\documentclass[11pt]{article}
\usepackage{amssymb}
\usepackage{amscd}
\input amssym.def
\input amssym  

\newtheorem{theorem}{Theorem}[section]
\newtheorem{corollary}{Corollary}[section]
\newtheorem{lemma}{Lemma}[section]

\newtheorem{remark}{Remark}[section]
\newtheorem{definition}{Definition}[section]
\newtheorem{proposition}{Proposition}[section]
\newtheorem{example}{Example}[section]

\usepackage{amssymb} 
\begin{document}

\title{Intertwining vertex operators and certain representations of
$\widehat{\mathfrak{sl}(n)}$ \footnote{ 2000 Mathematics Subject
Classification: 17B69, 17B65. Keywords: intertwining vertex operators,
standard modules, principal subspaces, recursions. }}

\author{ Corina Calinescu\footnote{The author gratefully acknowledges
partial support {}from NSF grant DMS 0100495 and from the Center for
Discrete Mathematics and Theoretical Computer Science (DIMACS),
Rutgers University.}}  \date{}

\maketitle

\renewcommand{\theequation}{\thesection.\arabic{equation}}
\renewcommand{\thetheorem}{\thesection.\arabic{theorem}}
\setcounter{equation}{0}
\setcounter{theorem}{0}
\setcounter{section}{0}

 \begin{abstract} We study the principal subspaces, introduced by
Feigin and Stoyanovsky, of the level $1$ standard modules for
$\widehat{\goth{sl}(l+1)}$ with $l \geq 2$. In this paper we construct
exact sequences which give us a complete set of recursions that
characterize the graded dimensions of the principal subspaces of these
representations. This problem can be viewed as a continuation of a new
program to obtain Rogers-Ramanujan-type recursions, which was
initiated by Capparelli, Lepowsky and Milas.  In order to
prove the exactness of the sequences we use intertwining vertex
operators and we supply a proof of the completeness of a list of
relations for the principal subspaces. By solving these recursions we
recover the graded dimensions of the principal subspaces, previously
obtained by Georgiev using a different method.  \end{abstract}

\section{Introduction}

The theory of vertex operator algebras (\cite{B}, \cite{FLM}, \cite
{FHL}; cf.\thinspace \cite{LL}) is a rich subject which provides,
among many other things, constructions of affine Lie algebras. In the
first vertex operator constructions, affine Lie algebras of type $A$,
$D$ and $E$, twisted by the principal automorphism of the underlying
finite-dimensional simple Lie algebra, were realized by means of
certain ``twisted vertex operator'' (\cite{LW1}, \cite{KKLW}; see also
\cite{L}). In \cite{FK} and \cite{S}, the untwisted affine algebras of
types $A^{(1)}$, $D^{(1)}$ and $E^{(1)}$ were represented using the
``untwisted vertex operators'' of dual-string theory. This point of
view was generalized and worked out in detail for all the standard
$\widehat{\goth{sl}(2)}$-modules in \cite{LW2}-\cite{LW4}, where the
${\cal Z}$-algebra program was initiated. As a consequence of the
${\cal Z}$-algebra theory, a proof of the classical Rogers-Ramanujan
identities was given and the Gordon-Andrews-Bressoud identities were
interpreted. These last identities were proved later in this spirit in
\cite{MP1}. Much work has been done in these directions (see
e.g. \cite{LP2}, \cite{Ca} and \cite{MP2}).

The sum sides of the two classical Rogers-Ramanujan identities are
  specializations of the solution of the Rogers-Ramanujan recursion
  (see the treatment in Chapter 7 of \cite{A}). By studying certain
  substructures, called the principal subspaces, of the level $1$
  standard $\widehat{\goth{sl}(2)}$-modules, Capparelli, Lepowsky
  and Milas recovered in \cite{CLM1} the Rogers-Ramanujan recursion
  and the graded dimensions (the generating functions of the
  dimensions of the homogeneous subspaces) of these subspaces.  In
  \cite{CLM2} they obtained the Rogers-Selberg recursions and showed
  that the graded dimensions of the principal subspaces of the level
  $k>1$ standard $\widehat{\goth{sl}(2)}$-modules satisfy these
  recursions.  In this paper we extend their approach to the principal
  subspaces of the level $1$ standard modules for
  $\widehat{\goth{sl}(l+1)}$, where $l \geq 2$. An announcement of
  some of the results proved in this paper is presented in
  \cite{C2}. In \cite{C3} we raise and solve similar problems for the
  principal subspaces of the level $k >1$ standard modules for
  $\widehat{\goth{sl}(3)}$.

The principal subspaces of the standard
  $\widehat{\goth{sl}(l+1)}$-modules were introduced and studied by
  Feigin and Stoyanovsky in \cite{FS1} and \cite{FS2}, and studied
  further by Georgiev in \cite{G1}.  These subspaces are generated by
  the affinization of the nilpotent subalgebra of $\goth{sl}(l+1)$
  consisting of the strictly upper triangular matrices and denoted by
  $\goth{n}$.  Feigin and Stoyanovsky showed that the duals of the
  principal subspaces are described in terms of symmetric polynomial
  forms not vanishing on certain hyperplanes. In this way they
  obtained the graded dimensions of the principal subspaces of the
  vacuum standard representations for $\widehat{\goth{sl}(2)}$ and
  $\widehat{\goth{sl}(3)}$.  Also, in the $\widehat{\goth{sl}(2)}$
  case, they expressed these graded dimensions in a different way, by
  using a geometric argument, and this led to the Rogers-Ramanujan and
  Gordon identities.  On the other hand, Georgiev computed differently
  the graded dimensions of the principal subspaces of certain standard
  $\widehat{\goth{sl}(l+1)}$-modules. He constructed combinatorial
  bases which are given in terms of partitions involving a color $j$
  for each simple root $\alpha_j$ ($1 \leq j \leq l$) and a charge
  $s$, $1 \leq s \leq k$ ($k$ being the level of the module) such that
  the parts of the same color and charge satisfy certain difference
  conditions. In the same program we would also like to mention
  \cite{P}, where Primc introduces Feigin-Stoyanovsky-type
  principal subspaces and constructs bases of them. Also, in
  \cite{ARS} Ardonne, Kedem and Stone use the principal
  subspaces to give formulas for the $q$-characters of arbitrary
  standard $\widehat{\goth{sl}(l+1)}$-modules.

  We consider the definition of the principal subspaces $W(\Lambda_i)$
  {}from \cite{FS1} and \cite{FS2}, namely
  $W(\Lambda_i)=U(\bar{\goth{n}}) \cdot v_{\Lambda_i}$, where
  $\bar{\goth{n}}=\goth{n} \otimes \mathbb{C}[t, t^{-1}]$ and
  $\Lambda_i$ and $v_{\Lambda_i}$ are the highest weights and highest
  weight vectors of the level $1$ standard representations
  $L(\Lambda_i)$ of $\widehat{\goth{sl}(l+1)}$ for $i=0, \dots ,
  l$. Of course, this definition works for any highest weight
  module. The space $W(\Lambda_i)$ can be identified with the quotient
  space $U(\bar{\goth{n}})/I_{\Lambda_i}$, where $I_{\Lambda_i}$ is
  the annihilator of $v_{\Lambda_i}$ in $U(\bar{\goth{n}})$.  In
  \cite{FS2} a result was announced describing these ideals (at least
  for the principal subspace of the vacuum representation
  $L(\Lambda_0)$), thus giving a presentation of the principal
  subspaces. In this paper we denote by
  $\widetilde{U(\bar{\goth{n}})}$ the completion of the universal
  enveloping algebra $U(\bar{\goth{n}})$ in the sense of \cite{LW3} or
  \cite{MP1}.

We consider the formal infinte sums
$$
R_{t}^{[j]}=\sum_{m_1+m_2=t} x_{\alpha_j}(m_1)x_{\alpha_j}(m_2),
$$
where $t \in \mathbb{Z}$ and $j=1, \dots, l$. Here we view
$R_{t}^{[j]}$ as elements of $\widetilde{U(\bar{\goth{n}})}$. Denote
by ${\mathcal J}$ the two-sided ideal of
$\widetilde{U(\bar{\goth{n}})}$ generated by the elements
$R_{t}^{[j]}$ for all $t \in \mathbb{Z}$ and $j=1, \dots, l$. We also
consider $\bar{\goth{n}}^{+}$ the subalgebra $\goth{n} \otimes
\mathbb{C}[t]$ of $\widehat{\goth{sl}(l+1)}$.  In this paper we prove
the following result: \\ \\ {\bf Theorem 1} {\it The annihilator of
  the highest weight vector of $L(\Lambda_0)$ in
  $U(\bar{\goth{n}})$ is described as follows:
$$
I_{\Lambda_0} \equiv {\mathcal J} \; \; \; \mbox{modulo} \; \; \; 
\widetilde{U(\bar{\goth{n}})\bar{\goth{n}}^{+}}.
$$
Moreover, for the principal subspaces $W(\Lambda_j)$ we have:   
$$
I_{\Lambda_j}\equiv {\mathcal J} + U(\bar{\goth{n}})x_{\alpha_j}(-1)
\; \; \; \mbox{modulo} \; \; \; 
\widetilde{U(\bar{\goth{n}})\bar{\goth{n}}^{+}},
$$}
where $j=1, \dots , l$.
\vspace{1em}

This theorem describes the left ideals $I_{\Lambda_j}$ of
$U(\bar{\goth{n}})$, and thereby gives a presentation of the principal
subspaces $W(\Lambda_j)$ for $j=0, \dots, l$.

More recently, a different proof of Theorem 1 has been given in
\cite{CalLM4}, which is part of an ongoing research program. New ideas
have been used in \cite{CalLM1} and \cite{CalLM2} to give a
presentation of the principal subspaces associated to the standard
modules for $\widehat{\goth{sl}(2)}$, and extended in \cite{CalLM3}
and \cite{CalLM4}.

By using Theorem 1 and intertwining vertex operators in the sense of
  \cite{FHL} and \cite{DL}, in this paper we find exact sequences for
  the principal subspaces of the standard
  $\widehat{\goth{sl}(l+1)}$-modules: \\ \\ {\bf Theorem 2} {\it The
  following sequences of maps between principal subspaces are exact:
$$
0 \longrightarrow W( \Lambda_1) \stackrel{e^{\lambda^1}}
\longrightarrow W( \Lambda_0) \stackrel{{\cal Y}_c (e^{\lambda_1}, x)} 
\longrightarrow W( \Lambda_1)
\longrightarrow 0,
$$
\[
0 \longrightarrow W(\Lambda_2) \stackrel{e^{\lambda^2}}
 \longrightarrow W( \Lambda_0) \stackrel{{\cal Y}_c (e^{\lambda_2},
 x)} \longrightarrow W( \Lambda_2) \longrightarrow 0, \]
 $$
 \vdots
 $$
  \[ 0 \longrightarrow W( \Lambda_l) \stackrel{e^{\lambda^l}}
 \longrightarrow W( \Lambda_0) \stackrel{{\cal Y}_c (e^{\lambda_l},
 x)} \longrightarrow W( \Lambda_l) \longrightarrow 0.  \] } \\ 
We
 refer to Section 2 for the linear maps ${\cal Y}_c(e_{\lambda_j}, x)$
 corresponding to intertwining operators. See Section 4 for the maps
 $e^{\lambda^j}$, where $\lambda^j$ are certain weights.

Let us denote by $\chi_{0}(x_1, \dots, x_l;q)$ the graded dimension of
   $W(\Lambda_0)$ with respect to certain natural weight and charge
   gradings, which are compatible.  The main consequence of the
   exactness of the sequences is: \\ \\ {\bf Theorem 3} {\it The
   graded dimension of the principal subspace $W(\Lambda_0)$ satisfies
   the following recursions:
   $$
   \chi_{0}(x_1, \dots, x_l;q)=\chi_{0}(x_1q, x_2, x_3, \dots , x_l;
   q)+x_1q\chi_{0}(x_1q^2, x_2q^{-1}, x_3, \dots, x_{l-1}, x_l;q),
   $$
   $$
   \chi_{0}(x_1, \dots, x_l; q)= \chi_{0}(x_1, x_2q, x_3, \dots , x_l;
   q)+x_2q \chi_{0}(x_1q^{-1}, x_2q^2, x_3q^{-1}, \dots, x_l;q),
   $$ 
   $$
   \vdots
   $$
   $$
   \chi_0(x_1, \dots, x_l;q)= \chi_0(x_1, x_2, x_3 , \dots, x_lq;q) +
   x_lq \chi_0 (x_1, x_2, x_3, \dots, x_{l-1}q^{-1}, x_lq^2;q).
   $$}

  These recursions can be reformulated in terms of the Cartan matrix
    as follows: \begin{eqnarray} \nonumber \chi_{0}(x_1, \dots, x_l;q)
    & = & \chi_{0}(x_1, \dots , (x_jq)^{\frac{a_{jj}}{2}}, \dots, x_l;
    q) \nonumber \\ &+
    &(x_jq)^{\frac{a_{jj}}{2}}\chi_{0}(x_1q^{a_{j1}}, x_2q^{a_{j2}},
    x_3q^{a_{j3}}, \dots, x_lq^{a_{jl}};q), \nonumber \end{eqnarray}
    where $A=(a_{ij})_{1 \leq i, j \leq l}$ is the Cartan matrix of
    $\goth{sl}(l+1)$.

Solving the recursions {}from Theorem 3 we obtain the graded dimension
of $W(\Lambda_0)$, and consequently the graded dimensions of all the
principal subspaces of the level $1$ standard modules for
$\widehat{\goth{sl}(l+1)}$.

The main feature of our strategy and results is understanding the
nature of the vertex-algebraic structure involved in the construction
of the exact sequences which give recursions.

We would like to mention that certain recurrence relations for the
 characters of the standard modules $L(k \Lambda_0)$ of any $A$, $D$,
 $E$ type Lie algebra have been obtained in \cite{CLiM}.

Although in this paper we concentrate on the case
$\widehat{\goth{sl}(l+1)}$, the generalization to affine algebras of
type $D$ and $E$ is possible and this will be the subject of a future
publication.

The paper is organized as follows. In Section 2 we recall the vertex
  operator construction of the level $1$ standard
  $\widehat{\goth{sl}(l+1)}$-modules and the construction of
  distinguished intertwining operators associated with the standard
  modules. Section 3 provides a proof of a presentation of the
  principal subspaces of these modules. In Section 4, the main section
  of this paper, we construct exact sequences using intertwining
  operators, find Rogers-Ramanujan-type recursions satisfied by the
  graded dimension of $W(\Lambda_0)$ and obtain the graded dimensions
  of $W(\Lambda_0), W(\Lambda_1), \dots, W(\Lambda_l)$.

This paper is part of the author's Ph.D. dissertation written under
the direction of James Lepowsky at Rutgers University.

{\bf Acknowledgement.} I want to acknowledge my deepest gratitude to
James Lepowsky for introducing me to this area of research, for his
enthusiasm and his inspiring guidance. I thank James Lepowsky and
Antun Milas for suggesting many improvements upon earlier versions of
this paper.

\section{Vertex operator constructions associated with standard
  modules}
\setcounter{equation}{0} 
In this section we recall the vertex operator
construction of the level $1$ standard
$\widehat{\goth{sl}(l+1)}$-modules for $l \geq 1$
(\cite{FK} and \cite{S}). We also review the vertex operator algebra
and module structures (\cite{B}, \cite{FLM}) associated to the level
$1$ standard modules and the construction of certain intertwining
operators among these modules (\cite{DL} and \cite{FHL}). We work in
the setting of \cite{FLM} (cf. \cite{LL}).

Let $\goth{g}$ be the finite-dimensional complex Lie algebra
$\goth{sl}(l+1)$ and let $\goth{h}$ be the Cartan subalgebra of
diagonal matrices.  We denote by $\Pi= \{ \alpha_1, \dots, \alpha_l
\}$ the set of simple (positive) roots, and by $\Delta_{+}$
(respectively, $\Delta_{-}$) the set of positive (respectively,
negative) roots of $\goth{g}$. Set $\Delta= \Delta_{+} \cup
\Delta_{-}$.  Denote by $\langle \cdot , \cdot \rangle$ the normalized
standard symmetric invariant nondegenerate bilinear form such that
$\langle \alpha, \alpha \rangle=2$ for any simple root $\alpha$ of
$\goth{g}$. For each root $\alpha$ fix a root vector $x_{\alpha}$.  We
fix $\{ x_{\alpha} \}_{\alpha \in \Delta} \cup \{h_{\alpha_i}
\}_{i=1}^{l}$ a Chevalley basis of $\goth{g}$.  Under our
identification we have $h_{\alpha_i}= \alpha_i$ for $i=1, \dots,
l$. The subalgebra $\sum_{i=1}^l \mathbb{C}h_{\alpha_i}$ of $\goth{g}$
is the Cartan subalgebra $\goth{h}$ mentioned before.  We identify
$\goth{h} \simeq \goth{h}^{*}$ via $\langle \cdot, \cdot \rangle$.

Let $Q=\sum_{i=1}^l \mathbb{Z} \alpha_i$ and $P=\sum_{i=1}^l
\mathbb{Z} \lambda_i$ be the root and weight lattices, where
$\lambda_i$ are the fundamental weights of $\goth{g}$ ($\langle
\lambda_i, \alpha_j \rangle = \delta_{i,j}$ for any $i,j =1, \dots,
l$).

It is known that there is a central extension of $P$ (and by
restriction this gives a central extension of $Q$) by a finite
(cyclic) group of roots of unity, which we denote by $A$, satisfying
the following condition:
$$
\epsilon (\alpha, \beta) \epsilon (\beta, \alpha)^{-1} = (-1)^{\langle
\alpha, \beta \rangle} \; \; \mbox{for} \; \; \alpha, \beta \in Q,
$$
where 
\begin{equation} \label{cocycle}
\epsilon: P \times P \longrightarrow A
\end{equation}
is a $2$-cocycle corresponding to the extension. Set
\begin{equation} \label{commutator} 
c(\lambda, \mu)= \epsilon (\lambda, \mu) \epsilon (\mu, \lambda)^{-1}
\; \; \mbox{for} \; \; \lambda, \mu \in P.
\end{equation}
This is the commutator map of the central extension (cf. \cite{FLM},
\cite{DL} and \cite{LL}).

We shall often be working with the positive nilpotent subalgebra of
$\goth{g}$,
\begin{equation} \label{nilpotent}
\goth{n}= \sum_{\alpha \in \Delta_{+}}\mathbb{C} x_{\alpha},
\end{equation}
which can be viewed as the subalgebra of $\goth{g}$ consisting of the
strictly upper triangular matrices.

Now we consider the untwisted affine Lie algebra associated to
$\goth{g}$,
\begin{equation}
\widehat{\goth{g}}= \goth{g} \otimes \mathbb{C}[t, t^{-1}] \oplus
\mathbb{C}c,
\end{equation}
where $c$ is a nonzero central element and 
\begin{equation}
[ x \otimes t^m, y \otimes t^n ] = [x, y] \otimes t^{m+n} + m\langle
x, y \rangle \delta _{m+n, 0} c
\end{equation}
for any $x, y \in \goth{g}$ and $m, n \in \mathbb{Z}$.  By adjoining
the degree operator $d$ ($[d, x \otimes t^m]=m$, $[d,c]=0$) to the Lie
algebra $\widehat{\goth{g}}$ one obtains the affine Kac-Moody algebra
$\widetilde{\goth{g}}=\widehat{\goth{g}} \oplus \mathbb{C}d$
(cf. \cite{K}).

Let us introduce the following subalgebras of $\widehat{\goth{g}}$:
\begin{equation} \label{n}
 \bar{\goth{n}}= \goth{n} \otimes \mathbb{C}[t, t^{-1}] ,
 \end{equation}
\begin{equation} \label{n+}
\bar{\goth{n}}^{+}=\goth{n} \otimes \mathbb{C}[t],
\end{equation}
\begin{equation}
\widehat{\goth{h}} = \goth{h} \otimes \mathbb{C}[t, t^{-1}] \oplus
\mathbb{C}c
\end{equation}
and
\begin{equation}
 \widehat {\goth{h}}_{\mathbb{Z}} = \coprod _{m \in \mathbb{ Z}
 \setminus {0}} \goth{h} \otimes t^m \oplus \mathbb{C}c.
 \end{equation} The latter is a Heisenberg subalgebra of $
 \widehat{\goth{g}}$ in the sense that its commutator subalgebra is
 equal to its center, which is one-dimensional. The form $\langle
 \cdot, \cdot \rangle $ on $\goth{h}$ extends naturally to $\goth{h}
 \oplus \mathbb{C}c \oplus \mathbb{C}d$. We shall identify $\goth{h}
 \oplus \mathbb{C}c \oplus \mathbb{C}d$ with its dual space $(\goth{h}
 \oplus \mathbb{C}c \oplus \mathbb{C}d)^{*}$ via this form.  The
 simple roots of $\widehat{\goth{g}}$ are $\alpha _{0}, \alpha_1,
 \dots, \alpha_l \in (\goth{h} \oplus \mathbb{C}c \oplus
 \mathbb{C}d)^{*}$.  We define $\Lambda_0, \dots, \Lambda_l \in (
 \goth{h} \oplus \mathbb{C}c \oplus \mathbb{C}d)^{*}$ by: $ \langle
 \Lambda_{i}, \alpha _{j} \rangle= \delta _{i, j}$, $ \langle
 \Lambda_{i}, c \rangle = 1$ for $i, j =0, \dots, l$, and $\langle
 \Lambda_0, d\rangle =0$, $\langle \Lambda_i, d \rangle= - 1/2 \langle
 \Lambda_i, \Lambda_i \rangle$ for $i = 1, \dots , l$. These are the
 fundamental weights of $\widehat{\goth{g}}$. Set $L( \Lambda_0),
 \dots , L(\Lambda_l)$ the standard $ \widehat {\goth{g}}$-modules of
 level 1 with highest weight vectors $v_{\Lambda_0}, \dots,
 v_{\Lambda_l}$. These modules are also called the integrable highest
 weight $\widehat{\goth{g}}$-modules of level $1$ (cf. \cite{K}).

{}From now on we will write $x(m)$ for the action of $x \otimes t^m$ on
any $\widehat{\goth{g}}$-module, where $x \in \goth{g}$ and $m \in
\mathbb{Z}$. Sometimes we will write $x(m)$ simply for the Lie algebra
element $x \otimes t^m$. It will be clear {}from the context whether
$x(m)$ is an operator or an element of $\widehat{\goth{g}}$.

We form the induced $ \widehat{\goth{h}}$-module
$$ M(1)= U( \widehat{\goth{h}}) \otimes _{U(\goth{h} \otimes
\mathbb{C}[t] \oplus \mathbb{C}c)}\mathbb{C},$$ such that $\goth{h}
\otimes \mathbb{C}[t]$ acts trivially and $c$ acts as identity on the
one-dimensional module $ \mathbb{C}$.  As a vector space, $M(1)$ is
isomorphic to the symmetric algebra $S(\widehat{\goth{h}}^{-})$, where
$\widehat{\goth{h}}^{-} = \goth{h} \otimes t^{-1}\mathbb{C} [t^{-1}]$.
The operators $\alpha(n)$, $n < 0$ act as mutiplication operators, the
operators $\alpha(n)$, $n>0$ act as certain derivations and $c$ acts
as identity.  Consider $\mathbb{C}[Q]$ and $\mathbb{C}[P]$ the group
algebras of the lattices $Q$ and $P$ with bases $\{ e^{\alpha} |
\alpha \in Q \}$ and $\{ e^{\lambda} | \lambda \in P \}$. Set
 $$ V_P = M(1) \otimes \mathbb{C}[P],$$
 $$V_Q= M(1) \otimes \mathbb{C}[Q],$$ and $$ V_Qe^{ \lambda_{j}}= M(1)
\otimes \mathbb{C}[Q]e^{ \lambda_{j}}, \; \; \; j=1, \dots, l.$$ Of
course, it is true that
 $$V_P= V_Q \oplus V_Qe^{ \lambda_{1}} \oplus \cdots \oplus
V_Qe^{\lambda_l}.$$

The vector space $V_{P}$ is equipped with certain structures. It has a
 structure of $\mathbb{C}[P]$-module, so that $e^{\lambda}$ with $
 \lambda \in P$ acts as the operator
 $$
 e^{\lambda}=1 \otimes e^{\lambda} \in \mbox{End} \;V_P.
 $$
By abuse of notation we shall often write $e^{\lambda}$ instead of $1
  \otimes e^{\lambda}$. In this paper, by $e^{\lambda}$ we mean either
  an operator or a vector of $V_P$, depending on the context.  Also,
  $V_P$ is an $\widehat{\goth{h}}$-module, such that the operators
  $\alpha(0)$ and $\alpha(n)$ for $n\neq 0$ act as follows:
 $$
 \alpha(0)(v \otimes e^{\lambda})= \langle \alpha,\lambda \rangle(v
 \otimes e^{\lambda}),
 $$
 $$\alpha(n)(v \otimes e^{\lambda})=\alpha(n) v \otimes e^{\lambda},
\; \; n\neq 0
 $$
 for $\alpha \in \goth{h} \simeq \goth{h}^{*}$ , $v \in M(1)$ and
 $\lambda \in P$.

Let $x, x_0, x_1, x_2 ,\dots$ be commuting formal variables.  For
 $\beta \in P$ we define $x^{\beta} \in (\mbox{End} \; V_P)\{\textit
 x\} $ (thought of as $x^{\beta(0)}$) by
 $$x^{\beta}(v \otimes e^{\lambda})=x^{ \langle \beta,\lambda
\rangle}(v \otimes e^{\lambda}),$$ where $v \in M(1)$ and $ \lambda
\in P$. By $(\rm{End}\; V_P) \{ \textit x \}$ we mean the space of all
formal series with complex (rather than integral) powers of $x$ and
with coefficients in $\rm{End} \; V_P$.

For every $\lambda \in P$ set 
\begin{equation} \label{str1}
Y(e^{\lambda},x)=E^{-}(-\lambda,x)E^{+}(-\lambda,x)
\epsilon_{\lambda}e^{\lambda}x^{\lambda},
\end{equation} an $\mbox{End} \; V_P$-valued formal Laurent series,
where $E^{\pm}(-\lambda, x)$ are defined as follows:
 $$
 E^{\pm}(-\lambda,x)=\mbox{exp}\left( \sum_{\pm n \geq
 1}\frac{-\lambda(n)}{n}x^{-n} \right )
 $$
 and
 $$
 \epsilon_{\lambda} e^{\beta}= \epsilon (\lambda, \beta) e^{\beta}
 $$
 for all $\alpha, \beta \in P$ (recall (\ref{cocycle})).  More
 generally, for a generic homogeneous vector $w \in V_P$ of the form
 $$
 w= \prod_{i=1}^{l}h_i(-n_i) \otimes e^{\lambda}, \; \; h_i \in
 \goth{h}, \; \; \lambda \in P \; \; \mbox{and} \; \; n_i \geq 1
 $$
 set \begin{equation} \label{operator} Y\left( \prod_{i=1}^l h_i(-n_i)
  \otimes e^{\lambda},x\right )= : \prod_{i=1}^l \left
  (\frac{1}{(n_i-1)} \left ( \frac{d}{dx} \right ) ^{n_i-1} h_i(x)
  \right )\; Y(e^{\lambda},x):.  \label{str2} \end{equation} Here
  $h_i(x)=\sum_{n \in \mathbb{Z}}h_i(n)x^{-n-1}$ and $ : \cdotp :$
  stands for a normal ordering operation.

 There is a natural $\widehat{\goth{sl}(l+1)}$-module structure on
 $V_P$ and this result is stated below.

 \begin{theorem} \label {FK} (\cite{FK}, \cite{S}; cf. \cite{FLM}) The
 vector space $V_P$ is an $ \widehat {\goth{sl}(l+1)}$-module of level
 $1$.  The action $x_{\alpha}(m)$ of $x_{\alpha} \otimes t^m$ for a
 root $\alpha$ and $m \in \mathbb{Z}$ is given by the coefficient of
 $x^{-m-1}$ in $ Y(e^{\alpha}, x)$. Moreover, the direct summands
 $V_Q$, $V_Qe^{ \lambda_{1}}, \dots , V_Qe^{ \lambda_{l}} $ of $V_P$
 are the standard $ \widehat {\goth{sl}(l+1)}$-modules of level $1$
 with highest weights $ \Lambda_{0}$, $\Lambda_{1}, \dots,\Lambda_{l}$
 and highest weight vectors $v_{\Lambda_0} =1 \otimes 1, v_{
 \Lambda_1} = 1 \otimes e^{ \lambda_{1}}, \dots , v_{\Lambda_l}= 1
 \otimes e^{ \lambda_{l}}$.
\end{theorem} 

We shall write 
\begin{equation} \label{vectors}
v_{\Lambda_0}= 1 \; \; \; \mbox{and} \; \; \;
v_{\Lambda_i}=e^{\lambda_i}
\end{equation}
for $i=1, \dots, l$.

For the sake of completeness we recall the following relations among
 the operators introduced above: \begin{proposition} (cf. \cite{LL})
 \label{comm} On $V_P$ we have \begin{equation}
 [\alpha(n),x^{\beta}]=0, \end{equation} \begin{equation}
 x^{\beta}e^{\lambda}=x^{\langle \beta,\lambda
 \rangle}e^{\lambda}x^{\beta}=e^{\lambda}x^{\beta+\langle
 \beta,\lambda \rangle}, \end{equation} \begin{equation}
 [\alpha(n),e^{\lambda}]=\delta_{n,0} \langle \alpha, \lambda \rangle
 e^{\lambda} \end{equation} and \begin{equation} \label{prod-oper}
 x_{\alpha}(m)e^{\lambda}=e^{\lambda}x_{\alpha}(m+\langle \alpha,
 \lambda \rangle )
\end{equation}
for any $\lambda, \beta \in P$, $\alpha \in \goth{h} \simeq
 \goth{h}^{*}$ and $m, n \in \mathbb{Z}$.  \end{proposition}

Now we recall the notion of vertex operator algebra, of module for
such a structure and of intertwining operator among given modules. We
refer to \cite {FLM}, \cite{FHL} and \cite{DL} (see also \cite {B},
\cite{LL}). In particular, a vertex operator algebra is a vector space
$V$ equipped with a vertex operator map 
\begin{eqnarray}
Y(\cdotp,x) : V & \longrightarrow & (\mbox{End V})[[x,x^{-1}]]
\nonumber \\ v & \mapsto & Y(v, x)=\sum_{n \in \mathbb{Z}} v_n
x^{-n-1}, \nonumber 
\end{eqnarray} a vacuum vector $1$ and a conformal vector $\omega$, such that 
the vector space $V$ has a $\mathbb{Z}$-grading truncated {}from below
with the homogeneous subspaces finite-dimensional, and such that the
axioms given in \cite{FLM} hold. The conformal vector $\omega$ gives
rise to the Virasoro algebra operators $L(n)$, $n \in \mathbb{Z}$:
$$Y(\omega,x)=\sum_{n \in \mathbb{Z}}L(n)x^{-n-2},$$ and the
$\mathbb{Z}$-grading of $V$ coincides with the eigenspace
decomposition given by the operator $L(0)$.

A module for the vertex operator algebra $V$ is a vector space
$W$ equipped with a vertex operator map 
\[ Y(\cdotp,x) : V
\longrightarrow (\mbox{End W})[[x,x^{-1}]], 
\] such that all the
axioms in the definition of the vertex operator algebra that make
sense hold, except that the grading on $W$ is allowed to be a
$\mathbb{Q}$-grading rather than a $\mathbb{Z}$-grading.

Suppose that $W_1$, $W_2$ and $W_3$ are modules for the vertex
operator algebra $V$. An intertwining operator of type
$$ \left( \begin{array} {c} W_1 \\ \begin{array}{cc} W_2 & W_3
\end{array} \end{array} \right) $$ is a linear map 
\begin{eqnarray}
\nonumber {\cal Y}(\cdotp, x) : W_2 & \longrightarrow & \mbox{Hom}
(W_3,W_1)\{x\} \nonumber \\ w &\mapsto & {\cal Y}(w,x)=\sum_{n\in
\mathbb{Q}}w_nx^{-n-1} \nonumber 
\end{eqnarray} that satisfies all the
defining properties of a module that make sense. The main axiom is the
Jacobi identity:

   \begin{eqnarray} \label{Jacobi-intertwining}
   \lefteqn{x_0^{-1}\delta \left ( \frac{x_1-x_2}{x_0} \right ) Y(u,
   x_1){\cal Y}(w_1, x_2)w_2} \nonumber \\
   &&\hspace{2em}-x_0^{-1}\delta \left ( \frac{x_2-x_1}{x_0} \right
   ){\cal Y}(w_1, x_2)Y(u, x_1)w_2 \nonumber \\
   &&\displaystyle{=x_2^{-1}\delta \left ( \frac{x_1-x_0}{x_2} \right
   ) {\cal Y}(Y(u, x_0)w_1, x_2)w_2} \nonumber \\ \end{eqnarray} for
   $u \in V$, $w_1 \in W_2$ and $w_2 \in W_3$.

  \begin{remark} \cite{FHL} \rm If the modules $W_1$, $W_2$ and $W_3$
  are irreducible (or indecomposable) then for any $w \in W_2$ we have
  ${\cal Y}(w, x) \in x^s \mbox{Hom}(W_3, W_1) \; [[x, x^{-1}]]$,
  where $s$ is a certain rational number. In this paper we will use
  intertwining operators among standard $\widehat{\goth{g}}$-modules
  for which $s$ will be an integer.  \end{remark}

 The standard modules $V_Q$ and $V_Qe^{\lambda_j}$ for $j=1, \dots, l$
 are endowed with natural vertex operator algebra and module
 structures.

\begin{theorem} (\cite{B}, \cite{FLM}) The formulas (\ref{str1}) and 
(\ref{str2}) give a vertex operator algebra structure on 
$V_Q$, and a $V_Q$-module structure on $V_P$. The vertex operator
algebra $V_Q$ is simple and $V_Qe^{ \lambda_{1}}, \dots , V_Q e^{
\lambda_{l}} $ are the irreducible modules for $V_Q$ (up to
equivalence).  \end{theorem}

{}From now on we will identify the standard modules as follows:
\begin{equation}
V_Q \simeq L(\Lambda_0), \; \; V_Q e^{ \lambda_{1}}
\simeq L(\Lambda_1), \dots,  V_Qe^{ \lambda_{l}} \simeq L(\Lambda_l).
\end{equation} 

We now discuss intertwining operators for the vertex operator algebra
 $L(\Lambda_0)$ and its standard modules $L(\Lambda_0), L(\Lambda_1),
 \dots, L(\Lambda_l)$. See Chapter 12 of \cite{DL} for a detailed
 treatment.  In order to construct intertwining operators one needs to
 define two operators $e^{i \pi \lambda}$ and $c(\cdotp, \lambda)$ on
 $V_P$ by:
 $$
e^{i \pi \lambda}(v \otimes e^{\beta})=e^{i \pi \langle \lambda ,
\beta \rangle}v \otimes e^{\beta},
 $$
 $$
 c(\cdotp, \lambda)(v \otimes e^{\beta})=c(\beta, \lambda)v \otimes
 e^{\beta},
 $$
 where $v \in M(1)$ and $ \beta , \lambda \in P$ (recall the
 commutator map (\ref{commutator})).  Let $r, s, p \in \{ 0, \dots,
 l\}$. Then \begin{eqnarray} { \cal Y}( \cdot , x): L(\Lambda_r) &
 \longrightarrow & \mbox{Hom}(L(\Lambda_s), L( \Lambda_p))\{ x \} \\ w
 & \mapsto & {\cal Y}(w, x)=Y(w, x)e^{i\pi\lambda_r}c(\cdot,
 \lambda_r) \nonumber \end{eqnarray} defines an intertwining operator
 of type \begin{equation} \label{type} \left( \begin{array} {c} L(
 \Lambda_p) \\ \begin{array}{cc} L( \Lambda_r) & L( \Lambda_s)
 \end{array} \end{array} \right).  \end{equation} We denote by ${\cal
 V}_{\Lambda_r, \Lambda_s}^{\Lambda_p}$ the vector space of the
 intertwining operators of type (\ref{type}). Set
 $$
 N_{\Lambda_r, \Lambda_s}^{\Lambda_p} = \mbox{dim} \; {\cal
 V}_{\Lambda_r, \Lambda_s}^{\Lambda_p}.
 $$
 These numbers are called {\it fusion rules} or {\it fusion
 coefficients}. It is known that $N_{\Lambda_r, \Lambda_s}^{\Lambda_p}
 \neq 0$ and in fact $N_{\Lambda_r, \Lambda_s}^{\Lambda_p}=1$, if and
 only if $p \equiv r+s \; \mbox{mod}\; (l+1)$ (cf. \cite{DL}).  The
 fact that all fusion rules in this situation are either zero or one
 follows also by using \cite{Li} since the modules $L(\Lambda_0),
 \dots, L(\Lambda_l)$ are simple currents.

In their work, Capparelli, Lepowsky and Milas
(\cite{CLM1}-\cite{CLM2}), and Georgiev (\cite{G1}) used certain
linear maps associated with intertwining operators. Similar maps will
play an important role in our paper.

For any $j=1, \dots, l$ we consider nonzero intertwining operators 
$$
{\cal Y}(\cdot , x): L(\Lambda_j) \longrightarrow \mbox{Hom} \;
(L(\Lambda_0), L(\Lambda_j)) [[x, x^{-1}]]
$$
(note that $N_{\Lambda_j, \Lambda_0} ^{\Lambda_j} =1$ cf \cite{DL}).
Also, ${\cal Y}( \cdot, x)$ involves only integral powers of $x$
(cf. Remark 5.4.2 of \cite{FHL}). In particular, we have a map
\begin{equation} \label{in}
{\cal Y}(e^{\lambda_j}, x): L(\Lambda_0) \longrightarrow
L(\Lambda_j)[[x, x^{-1}]].
\end{equation}

We denote the constant term (the coefficient of $x^0$) of ${\cal
Y}(e^{\lambda_j}, x)$ by ${\cal Y}_c(e^{\lambda_j}, x)$. This map
sends $v_{\Lambda_0}=1$ to a nonzero multiple of
$v_{\Lambda_j}=e^{\lambda_j}$. In this paper we normalize ${\cal
Y}(e^{\lambda_j}, x)$, if necessary, so that
\begin{equation} \label{int1}
{\cal Y}_c(e^{\lambda_j}, x) v_{\Lambda_0}= v_{\Lambda_j}.
\end{equation}

If we take $u=e^{\alpha}$ (for any positive root $\alpha$) and
$w_1=e^{\lambda_j}$ (for $j=1, \dots, l$) in the Jacobi identity
(\ref{Jacobi-intertwining}) and then apply $\mbox{Res}_{x_0}$ we
obtain
\begin{equation} 
[Y(e^{\alpha}, x_1), {\cal Y}(e^{\lambda_j}, x_2)]=0,
\end{equation}
which means that each coefficient of the series ${\cal
Y}(e^{\Lambda_j}, x)$ commutes with the action of $\bar{\goth{n}}$.
In particular, this is true for the constant term ${\cal
Y}_c(e^{\lambda_j}, x)$:
\begin{equation} \label{int2}
[Y(e^{\alpha}, x), {\cal Y}_c(e^{\lambda_j}, x)]=0.
\end{equation}
See \cite{DL}, \cite{G1}, \cite{CLM1} and \cite{CLM2} for further
discussion.

\section{Principal subspaces} 
\setcounter{equation}{0} 
The principal subspaces of the standard
modules $L(\Lambda_i)$, $i=0, \dots, l$ for $\widehat{\goth{sl}(l+1)}$
were introduced and studied by Feigin and Stoyanovsky in
\cite{FS1}-\cite{FS2}, and later by Georgiev in \cite{G1}. They are
generated by the affinization of the subalgebra $\goth{n}$ of
$\widehat{\goth{sl}(l+1)}$ consisting of the strictly upper triangular
matrices (recall (\ref{nilpotent})). These principal subspaces,
denoted by $W(\Lambda_i)$, are defined as follows: \begin{equation}
W(\Lambda_i) = U(\bar{\goth{n}}) \cdot v_{\Lambda_i}, \end{equation}
where $v_{\Lambda_i}$ are highest weight vectors of the standard
modules $L(\Lambda_i)$. By (\ref{vectors}) we also have
\begin{equation} \label{ps} W(\Lambda_i) = U(\bar{\goth{n}}) \cdot
e^{\lambda_i}.  \end{equation}

For $i=0, \dots, l$ let us consider the natural surjective maps
    \begin{eqnarray} \label{f maps} f_{\Lambda_i} : U(\bar{\goth{n}})
    & \longrightarrow & W(\Lambda_i) \\ a & \mapsto & a \cdot
    v_{\Lambda_i}. \nonumber \end{eqnarray} We denote by
    $I_{\Lambda_i}$ the annihilators of the highest weight vectors
    $v_{\Lambda_i}$ in $U(\bar{\goth{n}})$: \begin{equation}
    \label{ker} I_{\Lambda_i}=\mbox{Ker} \; f_{\Lambda_i}.
    \end{equation} These are left ideals in the associative algebra
    $U(\bar{\goth{n}})$. The aim of this section is to give a precise
    description of the ideals $I_{\Lambda_i}$. This is equivalent with
    finding a complete set of relations for $W(\Lambda_i)$, and thus
    with a presentation of these principal subspaces. This question
    was raised and discussed partially in \cite{FS2}. The presentation
    of the principal subspaces will play a key role in the process of
    obtaining recursions for the graded dimensions of $W(\Lambda_i)$.

    Denote by $\widetilde{U(\bar{\goth{n}})}$ the completion of
    $U(\bar{\goth{n}})$ in the sense of \cite{LW3} or \cite{MP1}.
    Recall {}from Theorem \ref{FK} the operators $x_{\alpha_j}(m)$,
    the images of $x_{\alpha_j} \otimes t^m$, for $j=1, \dots , l$ and
    $m \in \mathbb{Z}$. We consider the following formal infinite
    sums: \begin{equation} \label{sums}
      R_t^{[j]}=\sum_{m_1+m_2=t}x_{\alpha_j}(m_1)x_{\alpha_j}(m_2)
    \end{equation} for $t \in \mathbb{Z}$ and $j=1, \dots, l$,
    which act naturally on any $\widehat{\goth{g}}$-module. Denote by 
    ${\mathcal J}$ the two-sided ideal of $ \widetilde{U(\bar{\goth{n}})}$
 generated by 
    $R^{[j]}_t$ for $t \in \mathbb{Z}$ and $j=1, \dots , l$. Let $m$
       be an integer, possibly positive. To prove the presentation
       result it will be convenient to take $m_1, m_2 \leq m$, so that
       we truncate (\ref{sums}) as follows: \begin{equation}
       \label{R_t^j} R_{t;m}^{[j]}=\sum_{\begin{array} {c} m_1, m_2
       \leq m, \\
       m_1+m_2=t\end{array}}x_{\alpha_j}(m_1)x_{\alpha_j}(m_2).
       \end{equation} We shall often view $R_{t;m}^{[j]}$ as elements
       of $U(\bar{\goth{n}})$, rather than as endomorphisms of a
       $\widehat{\goth{g}}$-module. In general, it will be clear
       {}from the context when expressions such as (\ref{R_t^j}) are
       understood as elements of a universal enveloping algebra or as
       operators.  

       It is well known that the vector spaces $L(\Lambda_i)$ are
       graded with respect to a natural action of the Virasoro algebra
       operator $L(0)$, usually referred to as grading by {\it weight}
       (cf. \cite{LL}). For any integer $m$ and any root $\alpha$,
\begin{equation} 
\mbox{wt} \; x_{\alpha}(m)= -m,
\end{equation}
where $x_{\alpha}(m)$ is viewed as either an operator or as an element
of $U(\bar{\goth{n}})$. In particular, $U(\bar{\goth{n}})$ is graded
by weight. For any $\lambda \in P$ we have
\begin{equation}
\mbox{wt} \; e^{\lambda}= \frac{1}{2} \langle \lambda, \lambda \rangle.
\end{equation}
In particular,
\begin{equation}
\mbox{wt} \; v_{\Lambda_0}=0 \; \; \; \mbox{and} \; \; \; \mbox{wt} \;
v_{\Lambda_j}= \frac{1}{2} \langle\lambda_j, \lambda_j \rangle
\end{equation}
for all $j=1, \dots,l$ (recall (\ref{vectors})). See \cite{LL} for
further details and background.

We will restrict the weight gradings to the principal subspaces
$W(\Lambda_i)$. We have
\begin{equation} \label{direct sum}
W(\Lambda_0) = \coprod_{k \in \mathbb{Z}} W(\Lambda_0)_k 
\end{equation}
and
\begin{equation} \label {direct sum j}
W(\Lambda_j)= \coprod_{k \in \mathbb{Z}} W(\Lambda_j)_{k+
\frac{1}{2}\langle \lambda_j, \lambda_j \rangle},
\end{equation}
where $j=1, \dots , l$ and $W(\Lambda_0)_k$ and
$W(\Lambda_j)_{k+\frac{1}{2}\langle \lambda_j, \lambda_j \rangle}$ are
the subspaces generated by the homogeneous elements of weight $k$ and
$k+ \mbox{wt} \; v_{\Lambda_j}$, respectively.

\begin{remark} \label{ideals}
\rm
We have
\begin{equation}
L(0) I_{\Lambda_i} \subset I_{\Lambda_i}
\end{equation}
for $i=0, \dots, l$. Also, for any $j=1, \dots , l$ and $t , m \in
\mathbb{Z}$, $R_{t;m}^{[j]}$ has weight $-t$:
\begin{equation}
L(0) R_{t;m}^{[j]} =-t R_{t;m}^{[j]}.
\end{equation}
Also
$$
L(0) R_{t}^{[j]} =-t R_{t}^{[j]}.
$$
\end{remark}

We now prove the completeness of the list of relations for the
principal subspaces in the $\widehat{\goth{sl}(3)}$-case. We restrict
to this case $l=2$ for notational reasons only. The proof of the
corresponding theorem in the general case of
$\widehat{\goth{sl}(l+1)}$ is essentially identical, and we discuss it
at the end of this section. From now on until further notice we take
$l=2$ and we analyze the principal subspaces $W(\Lambda_i)$ of the
standard $\widehat{\goth{sl}(3)}$-modules $L(\Lambda_i)$ for
$i=0,1,2$. It will be clear that the next results have obvious
generalizations in the case $l \geq 2$, and that their proofs imitate
the proofs of the results {}from the case $l=2$.

The next theorem describes the left ideal $I_{\Lambda_0}$ of
$U(\bar{\goth{n}})$ and it was initially stated in \cite{FS2}. It also
gives a precise description of the left ideals $I_{\Lambda_1}$ and
$I_{\Lambda_2}$ of $U(\bar{\goth{n}})$ (recall (\ref{ker}) for the
notation).

\begin{theorem} \label{presentation sl3} The annihilators of the
  highest weight vectors of $L(\Lambda_i)$ for $i=0,1,2$ in
  $U(\bar{\goth{n}})$ are described as follows: \begin{equation}
    \label{structure}
    I_{\Lambda_0} \equiv {\mathcal J} \; \; \; \mbox{modulo} \; \; \; 
\widetilde{U(\bar{\goth{n}})\bar{\goth{n}}^{+}},
\end{equation} \begin{equation} \label{lambda_1} I_{\Lambda_1} \equiv 
  {\mathcal J}+U(\bar{\goth{n}})x_{\alpha_1}(-1) 
\; \; \; \mbox{modulo} \; \; \; 
\widetilde{U(\bar{\goth{n}})\bar{\goth{n}}^{+}},
\end{equation}
and
\begin{equation} \label{lambda_2}
I_{\Lambda_2} \equiv {\mathcal J}+ U(\bar{\goth{n}})x_{\alpha_2}(-1)
\; \; \; \mbox{modulo} \; \; \; 
\widetilde{U(\bar{\goth{n}})\bar{\goth{n}}^{+}},
\end{equation}
where ${\mathcal J}$ is the two-sided ideal of
$\widetilde{U(\bar{\goth{n}})}$ generated by $R_{t}^{[1]}$ and $R_{t}^{[2]}$
for $t \in \mathbb{Z}$.
    \end{theorem}

  In order to prove this theorem we need certain results which will be
  proved first.

  We set $\frak{B}=\{ x_{\alpha_2}, x_{\alpha_1},
  x_{\alpha_1+\alpha_2} \}$, a basis of the Lie subalgebra $\goth{n}=
  \mathbb{C}x_{\alpha_1} \oplus \mathbb{C}x_{\alpha_2} \oplus
  \mathbb{C}x_{\alpha_1+\alpha_2}$. We order the elements of
  $\frak{B}$ as follows:
 $$
 x_{\alpha_2} \prec x_{\alpha_1} \prec x_{\alpha_1+\alpha_2}.
 $$
 Set $\bar{\frak{B}}= \{ x(n) \; | \; x \in {\frak{B}}, \; n \in
 \mathbb{Z} \}$, so that $\bar{\frak{B}}$ is a basis of the Lie
 algebra $\bar{\goth{n}}$. We choose the following total order
 $\preceq$ on $\bar{\frak{B}}$:
 $$
 x_1(i_1) \preceq x_2(i_2) \; \; \mbox{iff} \; \; x_1 \prec x_2 \; \;
 \mbox{or} \; \; x_1=x_2 \; \; \mbox{and} \; \; i_1 \leq i_2.
 $$
 Then by the Poincar\'e-Birkhoff-Witt theorem we obtain a basis of the
 universal enveloping algebra $U(\bar{\goth{n}})$:
  $$
x_{\alpha_2}(n_1) \cdots x_{\alpha_2}(n_s)x_{\alpha_1}(m_1) \cdots
x_{\alpha_1}(m_r)x_{\alpha_1+\alpha_2}(l_1) \cdots
x_{\alpha_1+\alpha_2}(l_p)
 $$
with $n_1 \leq \cdots \leq n_s$, $m_1 \leq \cdots \leq m_r$ and $l_1
\leq \cdots \leq l_p$.

We shall refer to the expressions $x_{\gamma_1}(m_1) \cdots
x_{\gamma_r}(m_r)$ with $m_1, \dots , m_r \in \mathbb{Z}$ and
$\gamma_1, \dots, \gamma_r \in \{ \alpha_1, \alpha_2,
\alpha_1+\alpha_2 \}$ as the (noncommutative) {\it monomials} in
$U(\bar{\goth{n}})$.

   \begin{remark} \label{homo} \rm Consider the monomial
   $$x_{\alpha_2}(n_1) \cdots x_{\alpha_2}(n_s)x_{\alpha_1}(m_1)
   \cdots x_{\alpha_1}(m_r)x_{\alpha_1+\alpha_2}(l_1) \cdots
   x_{\alpha_1+\alpha_2}(l_p)$$ with $m_r, l_p \leq -1$, a homogeneous
   (with respect to the weight grading) basis element of
   $U(\bar{\goth{n}})$. By using the Lie brackets this monomial can be
   written as a linear combination of monomials $x_{\alpha_2}(m_{1,2})
   \cdots x_{\alpha_2}(m_{r_{2},2})x_{\alpha_1}(m_{1,1}) \cdots
   x_{\alpha_1}(m_{r_{1},1})$ of the same weight, such that $m_{1,2}
   \leq \cdots \leq m_{r_{2},2}$, $m_{1,1} \leq \cdots \leq
   m_{r_{1},1} \leq -1$, and elements of the left ideal
   $U(\bar{\goth{n}})\bar{\goth{n}}^{+}$. This result was proved and
   exploited in \cite{G1}.  \end{remark}

\begin{remark} \label{negative weight}
\rm Note that any homogeneous element $a \in U(\bar{\goth{n}})$ of
nonpositive weight lies in
$U(\bar{\goth{n}})\bar{\goth{n}}^{+}$. Indeed, $a$ can be written as a
linear combination of monomials of the same weight, and by using
brackets if necessary, we see that each monomial lies in
$U(\bar{\goth{n}}) \bar{\goth{n}}^{+}$.
\end{remark}

Next we will show that only certain monomials of the form
$$x_{\alpha_2}(m_{1,2}) \cdots x_{\alpha_2}(m_{r_2,
2})x_{\alpha_1}(m_{1,1}) \cdots x_{\alpha_1}(m_{r_1,1}) \in
U(\bar{\goth{n}})$$ applied to the highest weight vector
$v_{\Lambda_i}$ are needed to span the principal subspace
$W(\Lambda_i)$, $i=0, 1, 2$. We will first do a straightforward
computation for any $m \in \mathbb{Z}$,
 $$
 x_{\alpha_2}(1)x_{\alpha_1}(m)-x_{\alpha_2}(0)
x_{\alpha_1}(m+1)=x_{\alpha_1}(m)x_{\alpha_2}(1)-
x_{\alpha_1}(m+1)x_{\alpha_2}(0)
\in U(\bar{\goth{n}})\bar{\goth{n}}^{+},
 $$
and this observation leads to the following:

\begin{lemma}  \label{truncation}
Let $r_1, r_2 \geq 1$ and let $m_{1,2}, \dots, m_{r_2,2}, m_{1,1},
\dots, m_{r_1,1}$ be integers such that $m_{1,2} \leq \cdots \leq
m_{r_2,2}$ and $m_{1,1} \leq \cdots \leq m_{r_1,1} \leq -1$. Assume
that $m_{r_2,2} \geq r_1$. Then
\begin{equation} \label{elem}
x_{\alpha_2}(m_{1,2}) \cdots x_{\alpha_2}(m_{r_2,2})
x_{\alpha_1}(m_{1,1}) \cdots x_{\alpha_1}(m_{r_1,1})
\end{equation}
is a linear combination of monomials of the same weight as (\ref{elem}),
\begin{equation}
x_{\alpha_2}(m_{1,2}') \cdots x_{\alpha_2}(m_{r_2,2}')
 x_{\alpha_1}(m_{1,1}') \cdots x_{\alpha_1}(m_{r_1,1}'),
 \end{equation} with \begin{equation} \label{cond0} m_{1,2}' \leq
 \cdots \leq m_{r_2,2}' \leq r_1-1, \; m_{1,1}' \leq \cdots \leq
 m_{r_1,1}' \leq -1, \end{equation} and monomials of
 $U(\bar{\goth{n}})\bar{\goth{n}}^{+}$. In particular, any homogeneous
 element of $U(\bar{\goth{n}})$ can be written as a sum of two
 homogeneous elements, one such that the sequences of integers satisfy
 (\ref{cond0}) and the other an element of
 $U(\bar{\goth{n}})\bar{\goth{n}}^{+}$.  \end{lemma} {\em Proof:} We
 may assume that the monomial (\ref{elem}) has positive
 weight. (Otherwise, it lies in $U(\bar{\goth{n}})\bar{\goth{n}}^{+}$
 by Remark \ref{negative weight}.)  We first show that for any
 integers $m_{r_2,2}, m_{1,1}, \dots , m_{r_1,1}$ such that $m_{r_2,1}
 \geq r_1$ and $m_{1,1} \leq \cdots \leq m_{r_1,1} \leq -1$ we have
 \begin{eqnarray} \nonumber &&
 x_{\alpha_2}(m_{r_2,2})x_{\alpha_1}(m_{1,1}) \cdots
 x_{\alpha_1}(m_{r_1,1}) \nonumber \\ &&- \sum_{t=1}^{r_1}
 x_{\alpha_2}(m_{r_2,2}-1) x_{\alpha_1}(m_{1,1}) \cdots
 x_{\alpha_1}(m_{t, 1}+1) \cdots x_{\alpha_1}(m_{r_1,1}) \nonumber \\
 &&+ \sum_{\stackrel{s,t=1}{r\neq s}}^{r_1} x_{\alpha_2}(m_{r_2,2}-2)
 x_{\alpha_1}(m_{1,1}) \cdots x_{\alpha_1}(m_{s,1}+1) \cdots
 x_{\alpha_1}(m_{t,1}+1) \cdots x_{\alpha_1}(m_{r_1,1})\nonumber \\ &&
 \hspace{15em} \vdots \nonumber \\ &&+(-1)^{r_1}
 x_{\alpha_2}(m_{r_2,2}-r_1) x_{\alpha_1}(m_{1,1}+1) \cdots
 x_{\alpha_1}(m_{r_1,1}+1) \in
 U(\bar{\goth{n}})\bar{\goth{n}}^{+}. \label{expr} \end{eqnarray}
 Indeed, we use the brackets to move $x_{\alpha_2}(m_{r_2,2}-i)$ for
 $0 \leq i \leq r_1$ to the right of their corresponding monomials. At
 each step we cancel the similar terms with opposite signs (note the
 alternation of the signs in (\ref{expr})). Thus the monomial
 $x_{\alpha_2}(m_{r_2,2}) x_{\alpha_1}(m_{1,1}) \cdots
 x_{\alpha_1}(m_{r_1,1})$ is a linear combination of monomials
 $$x_{\alpha_2}(m_{r_2,2}') x_{\alpha_1}(m_{1,1}') \cdots
 x_{\alpha_1}(m_{r_1,1}'),$$ where $m_{r_2,2}' \leq r_1-1$ and
 $m_{1,1}' \leq \cdots m_{r_1,1}' \leq -1$ and monomials of
 $U(\bar{\goth{n}})\bar{\goth{n}}^{+}$.  In particular, this
 proves the lemma when $r_2=1$.

Now assume that $r_2>1$.  We multiply each $x_{\alpha_2}(m_{r_2,2}')
x_{\alpha_1}(m_{1,1}') \cdots x_{\alpha_1}(m_{r_1,1}')$ to the left by
$x_{\alpha_2}(m_{1,2}) \cdots x_{\alpha_2}(m_{r_{2}-1,2})$ and arrange
the indices of the sequence $m_{1,2}, \dots, m_{r_2-1, 2}, m_{r_2,2}'$
in a nonincreasing order {}from right to left. If all indices of this
sequence are less than $r_1$ then we are done.  Assume that there is
at least one index which is greater than $r_1-1$. We apply the same
strategy as in the previous case and we obtain the desired conclusion
after a finite number of steps.  $\Box$
\vspace{1em}

We also have:

\begin{lemma} \label{truncation 12}
 Let $r_1, r_2 \geq 1$ and $m_{1,2}, \dots, m_{r_2,2}, m_{1,1}, \dots,
 m_{r_1,1}$ be integers such that $m_{1,2} \leq \cdots \leq m_{r_2,2}$
 and $m_{1,1} \leq \cdots \leq m_{r_1,1} \leq -1$, and consider
 \begin{equation} \label{el} x_{\alpha_2}(m_{1,2}) \cdots
 x_{\alpha_2}(m_{r_2,2}) x_{\alpha_1}(m_{1,1}) \cdots
 x_{\alpha_1}(m_{r_1,1}) \in U(\bar{\goth{n}}).
\end{equation} 
\begin{enumerate}
 \item Assume that $m_{r_2,2} \geq r_1$. Then (\ref{el})
is a linear combination of monomials of the same weight
\begin{equation} \label{ele}
x_{\alpha_2}(m_{1,2}') \cdots x_{\alpha_2}(m_{r_2,2}')
 x_{\alpha_1}(m_{1,1}') \cdots x_{\alpha_1}(m_{r_1,1}') \end{equation}
 with
 $$
 m_{1,2}' \leq \cdots \leq m_{r_2,1}' \leq r_1-1, \; m_{1,1}' \leq
 \cdots \leq m_{r_1,1}' \leq -2,
 $$
 and monomials of $U(\bar{\goth{n}})\bar{\goth{n}}^{+}$ and of 
$U(\bar{\goth{n}})x_{\alpha_1}(-1)$.

 \item Assume that $m_{r_2,2} \geq r_1-1$. Then (\ref{el}) is a linear
 combination of monomials as in (\ref{ele}) such that
 $$
 m_{1,2}' \leq \cdots \leq m_{r_2,2}' \leq r_1-2, \; m_{1,1} \leq
 \cdots \leq m_{r_1,1}' \leq -1,
 $$
 and monomials in $U(\bar{\goth{n}})\bar{\goth{n}}^{+}$ and in 
$U(\bar{\goth{n}})x_{\alpha_2}(-1)$.
\end{enumerate}
In particular, the assertions (1) and (2) hold for homogeneous
 elements. $\; \; \; \Box$ \end{lemma}

Since the proofs of Lemmas \ref{truncation} and \ref{truncation 12}
are similar, we omit the proof of Lemma \ref{truncation 12}.

 In view of the above results we see that for $i=0,1,2$ the principal
 subspaces $W(\Lambda_i)$ are spanned by
\begin{equation}
x_{\alpha_2}(m_{1,2}) \cdots
x_{\alpha_2}(m_{r_2,2})x_{\alpha_1}(m_{1,1}) \cdots
x_{\alpha_1}(m_{r_1,1}) \cdot v_{\Lambda_i},
\end{equation}
such that
\begin{equation}
m_{1,2} \leq \cdots \leq m_{r_2,2} \leq r_1-1- \delta_{i, 2}, \;
\mbox{and} \; m_{1,1} \leq \cdots \leq m_{r_1,1}\leq -1- \delta_{i,1}.
\end{equation}

Inspired by \cite{LW3} we introduce the following strict linear
        ordering. It is a lexicographical ordering on the cartesian
        product of two ordered sets of sequences of lengths $r$ and
        $s$.  \begin{definition} \label{lexicographic} \rm Let $r$ and
        $s$ be positive integers and $m_1, \dots , m_{r+s}$ integers
        such that \begin{equation} \label{c1} m_1 \leq \cdots \leq
        m_r, \; \; m_{r+1} \leq \cdots \leq m_{r+s} \end{equation} and
        \begin{equation} \label{c2} m_1+ \cdots + m_r+m_{r+1} + \cdots
        +m_{r+s}=k \end{equation} for a fixed integer $k$.  Let $m_1',
        \dots , m_r', m_{r+1}', \dots , m_{r+s}'$ be another sequence
        of integers which satisfies the conditions (\ref{c1}) and
        (\ref{c2}).  We write
    $$
    (m_1', \dots , m_r', m_{r+1}', \dots , m_{r+s}') < (m_1, \dots,
    m_r, m_{r+1} , \dots, m_{r+s})
    $$
    if there exists $t \in \mathbb{Z}_{+}$ with $1 \leq t \leq r+s$
    such that $m_1'=m_1, \dots, m_{t-1}'=m_{t-1}$ and $m_t' < m_t$.
    \end{definition}

    If we require in the previous definition that $m_1 \leq \cdots
    \leq m_r \leq k'$ and $m_{r+1} \leq \cdots \leq m_{r+s} \leq k''$
    for some fixed integers $k'$ and $k''$ then there are only
    finitely many sequences which are less than the initial one in the
    lexicographical ordering.  We say that the set of monomials
    $$
    \{ x_{\alpha_2}(m_1) \cdots x_{\alpha_2}(m_r)x_{\alpha_1}(m_{r+1})
    \cdots x_{\alpha_1}(m_{r+s}) \}
    $$
   with
    $$
m_1 \leq \cdots \leq m_r, \; m_{r+1} \leq \cdots \leq m_{r+s}, \;
\mbox{and} \; m_1+ \cdots + m_{r+s}=k
    $$
    is linearly ordered by ``$<$'' if we use Definition
    \ref{lexicographic} for their corresponding sequences of indices.

Let $m_{1, 2}, \dots, m_{r_2, 2}, m_{1,1}, \dots m_{r_1, 1} \in
    \mathbb{Z}$. We say that the monomial
    $$x_{\alpha_2}(m_{1, 2}) \cdots x_{\alpha_2}(m_{r_2,
2})x_{\alpha_1} (m_{1, 1}) \cdots x_{\alpha_1}(m_{r_1, 1})$$ satisfies
{\it the difference two condition} if the sequences $m_{1, 2} \dots,
m_{r_2, 2}$, $m_{1,1}, \dots , m_{r_1, 1}$ satisfy the difference two
condition, i.e.,
    $$
    m_{s, 2} - m_{s-1, 2} \geq 2 \; \; \; \mbox{for any} \; \; s, \;
    \; 2 \leq s \leq r_2
    $$
  and
  $$
  m_{t, 1}- m_{t-1, 1} \geq 2 \; \; \; \mbox{for any} \; \; t, \; \; 2
  \leq t \leq r_1.
  $$

     \begin{lemma} \label{p1} Let $r_1, r_2 \geq 1$ and $m_{1,2},
    \dots , m_{r_2,2}, m_{1,1}, \dots , m_{r_1,1}$ be integers such
    that $m_{1,2} \leq \cdots \leq m_{r_2,1} \leq r_1-1$ and $m_{1,1}
    \leq \cdots \leq m_{r_1,1} \leq -1$. Then
$$
    M=x_{\alpha_2}(m_{1,2})\cdots
    x_{\alpha_2}(m_{r_2,2})x_{\alpha_1}(m_{1,1}) \cdots
    x_{\alpha_1}(m_{r_1,1})
    $$ can be expressed as
    $$
    M=M_1+M_2+M_3,
    $$
    where $M_1$ is a linear combination of monomials satisfying the
    difference two condition, $M_2 \in {\mathcal J}$, $M_3 \in
    \widetilde{U(\bar{\goth{n}})\bar{\goth{n}}^{+}}$, and where
    $\rm{wt} \; M_1= \rm{wt} \; M_2= \rm{wt} \; M_3= \rm{wt} \; M$. In
    particular, any homogeneous element $\Phi \in U(\bar{\goth{n}})$
    satisfying the above conditions for their corresponding sequences
    of integers is a linear combination of the elements $\Phi_1, \Phi_2, \Phi_3$,
    where $\Phi_1$ a linear combination of monomials that satisfy the
    difference two condition, $\Phi_2 \in {\mathcal J}$ and $\Phi_3
    \in \widetilde{U(\bar{\goth{n}})\bar{ \goth{n}}^{+}}$, and such
    that $\Phi$, $\Phi_1$, $\Phi_2$ and $\Phi_3$ have the same
    weight.  \end{lemma} {\em Proof:} Assume that $M$ has positive
  weight. If $\mbox{wt} \; M \leq 0$ then $M \in
  U(\bar{\goth{n}})\bar{\goth{n}}^{+}$ by Remark \ref{negative
    weight}. In order to prove this lemma we study four situations.

       Case I. Assume that there is at least one index in the sequence
       $m_{1,2} \cdots, m_{r_2,2}$, say $m_{s,2}$, such that $m_{s,2}
       +m_{1,1} \geq 0$. We move $x_{\alpha_2}(m_{s,2})$ to the right
       of $M$ by using brackets. Observe that $M$ decomposes as a sum
       of monomials which end in
       $x_{\alpha_1+\alpha_2}(m_{s,2}+m_{1,1}), \dots ,
       x_{\alpha_1+\alpha_2}(m_{s,2}+m_{r_1,1})$ and
       $x_{\alpha_2}(m_{s,2})$. Note that $m_{s, 2}+ m_{1,1}, \dots,
       m_{s,2}+m_{r_1, 1}, m_{s,2} \geq 0$. Hence $M \in
       U(\bar{\goth{n}})\bar{\goth{n}}^{+}$.

   Case II. Suppose that $r_2>r_1$ and that at least $r_1+1$ indices
   of the sequence $m_{1,2}, \dots, m_{r_2,2}$ are
   nonnegative. Without loss of generality we may and do assume that
   $m_{r_2-r_1, 2}, \dots, m_{r_2,2} \geq 0$. As in the previous case,
   $M$ can be written as a sum of monomials which end in
   $x_{\alpha_2}(m_{s,2})$ for $r_2-r_1 \leq s \leq r_2$. Since each
   $m_{s, 2} \geq 0$, thus $M \in U(\bar{\goth{n}})
   \bar{\goth{n}}^{+}$.

   Before we analyze the remaining situations we would like to show
   that any monomial of fixed positive weight $x_{\alpha}(m_1) \cdots
   x_{\alpha}(m_r)$, such that $m_1 \leq \cdots \leq m_r \leq -1$ and
   $\alpha= \alpha_1$ or ($\alpha=\alpha_2$) is a linear combination
   of elements in the ideal ${\mathcal J}$ generated by $R_{t}
   =R_{t}^{[1]}$ (or $R_{t}=R_{t}^{[2]}$) for $t<0$ modulo
   $\widetilde{U(\bar{\goth{n}})\bar{\goth{n}}^{+}}$, and monomials of
   the same weight which satisfy the difference two condition.

       We start with the first ({}from right to left) pair of
       consecutive indices $(m_{s-1}, m_s)$ of the sequence $(m_1,
       \dots, m_r)$ such that $m_s-m_{s-1} \leq 1$. We replace this
       pair by
      $$
       R_{m_{s-1}+m_s; -1}-\sum_{m_{s}'- m_{s-1}' \geq 2}
       x_{\alpha}(m_{s-1}')x_{\alpha}(m_s')
       $$  
       with $m_{s-1}'+m_s'=m_{s-1}+m_s$. Thus \[ x_{\alpha}(m_1)
       \cdots x_{\alpha}(m_r)- \sum_{\textbf m} A_{\textbf m}
       x_{\alpha}(m_1) \cdots x_{\alpha}(m_{s-1}') x_{\alpha}(m_s')
       \cdots x_{\alpha}(m_r) \in {\mathcal J} 
       \; \; \mbox{modulo} \; \; 
       \widetilde{U(\bar{\goth{n}})\bar{\goth{n}}^{+}},\]  where
       ${\textbf m}=(m_1, \dots, m_{s-1}', m_s', \dots, m_r)$ and
       $A_{\textbf m}$ are constants. We arrange the indices of
       ${\textbf m}$ in a nonincreasing order {}from right to left and
       we observe that each sequence ${\textbf m}$ is less than the
       initial one $(m_1, \dots, m_r)$ in the lexicographic order.
       Applying the same procedure, in a finite number of steps, we
       are able to write $x_{\alpha}(i_1) \cdots x_{\alpha}(i_r)$ as a
       linear combination of monomials which satisfy the difference
       two condition and elements in the ideal ${\mathcal J}$ modulo
       $\widetilde{U(\bar{\goth{n}})\bar{\goth{n}}^{+}}$.

      Case III. Suppose that $m_{1,2} \leq \cdots \leq m_{r_2,2} \leq
      -1$. We proceed as above and use $R_{t}^{[2]}$ to
      obtain
   $$  
  M- \sum_{\textbf m_2, \; \textbf m_1} A_{\textbf m_2, \; \textbf
  m_1} x_{\alpha_2}(m_{1,2}') \cdots
  x_{\alpha_2}(m_{r_2,2}')x_{\alpha_1}(m_{1,1}) \cdots
  x_{\alpha_1}(m_{r_1,1}) \in {\mathcal J} \; \; \mbox{modulo} \; \;  
  \widetilde{U(\bar{\goth{n}})\bar{\goth{n}}^{+}},
  $$
where $A_{\textbf m_2, \; \textbf m_1}$ are constants and
  $$
  {\textbf m_2}= (m_{1,2}', \dots, m_{r_2,2}'), \; \;
  m_{s,2}'-m_{s-1,2}' \geq 2\; \; \mbox{for} \; \; 2 \leq s \leq r_2,
  \; \; m_{r_2,1}' \leq -1
  $$
  and
  $$
  {\textbf m_1}= (m_{1,1}, \dots, m_{r_1, 1}), \; \; m_{1,1} \leq
  \cdots \leq m_{r_1, 1} \leq -1.
  $$
  Then we use the expressions $R_{t}^{[1]}$ with $t<0$ to make
difference two among the indices of the sequence ${\textbf m_1}$. The
conclusion follows.

Case IV. Now assume that there is at least $1$ and at most $r_1$
nonnegative integers in the sequence $m_{1,2},\dots, m_{r_2,2}$. In
this case it is essential to use first $R_{t}^{[2]}$ for $t \leq
2r_1-2$, to make difference two among $m_{1,2}, \dots, m_{r_2,2}$. If
our new sequences have indices $m_{s,2}$ such that $m_{s,2} \geq r_1$
we apply Lemma \ref{truncation} and we also apply Case I if $m_{s,2}
+m_{1,1} \geq 0$ for $1 \leq s \leq r_2$.  Hence, we have
  $$ 
  M- \sum_{\textbf m_2, \; \textbf m_1} A_{{\textbf m_2}, \; {\textbf
  m_1}} x_{\alpha_2}(m_{1,2}') \cdots
  x_{\alpha_2}(m_{r_2,2}')x_{\alpha_1}(m_{1,1}') \cdots
  x_{\alpha_1}(m_{r_1, 1}') \in {\mathcal J} + \widetilde{U(\bar{\goth{n}})
  \bar{\goth{n}}^{+}},
  $$                   
  where $A_{\textbf m_2, \; \textbf m_1}$ are constants and
  $${\textbf m_2}= (m_{1,2}' \dots, m_{r_2,2}'), \; \; \; { \textbf
m_1}= (m_{1,1}', \dots, m_{r_1,1}'),
  $$
   $$ m_{1,2}' \leq \cdots \leq m_{r_2,2}' \leq r_1-1, \; \; m_{s,2}'
-m_{s-1, 2}' \geq 2 \; \; \mbox{for} \; \; 2 \leq s \leq r_2,
   $$
   and
   $$
    m_{1,1}'\leq \cdots \leq m_{r_1,1}' \leq -1.
    $$
Notice that the sequences $({\textbf m_2}, \; {\textbf m_1})$ are less
     than $(m_{1,2}, \dots, m_{r_2,2}, m_{1,1}, \dots , m_{r_1,1})$ in
     the lexicographic ordering $``<"$. In order to finish the proof
     of this case we use the expressions $R_{t}^{[1]}$ for $t \leq
     -2$ to make difference two among the indices $m_{1,1}', \dots ,
     m_{r_1,1}'$. We finally observe that in all situations the
     weights of $M_1, M_2, M_3$ equal the weight of $M$.  $\: \: \:
     \Box$ \vspace{1em}

As an illustration of the technique used in the proof of the previous
   proposition we give a typical example.  \begin{example} \rm We show
   that the monomial
         $$
         M=x_{\alpha_2}(0)x_{\alpha_2}(1)x_{\alpha_1}(-3)x_{\alpha_1}(-3)
         $$ 
         of weight $5$ has the decomposition stated  in Lemma \ref{p1}.

         First we apply Case IV to write $M$ as follows:
$$M=\frac{1}{2}R^{[2]}_{1}x_{\alpha_1}(-3)x_{\alpha_1}(-3)-
x_{\alpha_2}(-1)x_{\alpha_2}(2)x_{\alpha_1}(-3)x_{\alpha_1}(-3)-N,
         $$
         with $N \in \widetilde{U(\bar{\goth{n}})\bar{\goth{n}}^{+}}$ by Case I.
         By using Lemma \ref{truncation} we have \begin{eqnarray}
           \nonumber \lefteqn{x_{\alpha_2}(-1)x_{\alpha_2}(2)
             x_{\alpha_1}(-3)x_{\alpha_1}(-3)} \nonumber \\ &&
           =2x_{\alpha_2}(-1)x_{\alpha_2}(1)x_{\alpha_1}(-3)x_{\alpha_1}(-2)
           \nonumber \\ &&\hspace{1em}+x_{\alpha_2}(-1)x_{\alpha_2}(0)
           x_{\alpha_1}(-2)x_{\alpha_1}(-2) +P,
           \nonumber \end{eqnarray} where $P \in
         U(\bar{\goth{n}})\bar{\goth{n}}^{+}$.  Finally, we observe
         that \begin{eqnarray} \nonumber
           \lefteqn{x_{\alpha_2}(-1)x_{\alpha_2}(1)x_{\alpha_1}(-3)
             x_{\alpha_1}(-2)} \nonumber \\
           &&=\frac{1}{2}x_{\alpha_2}(-1)x_{\alpha_2}(1)R^{[1]}_{-5} 
           \nonumber \\ && \hspace{1em}-\displaystyle
           x_{\alpha_2}(-1)x_{\alpha_2}(1)x_{\alpha_1}(-4)x_{\alpha_1}(-1)
           \nonumber 
         \end{eqnarray} 
and 
\begin{eqnarray} \nonumber
           &&x_{\alpha_2}(-1)x_{\alpha_2}(0)x_{\alpha_1}(-2)x_{\alpha_1}(-2)
           \nonumber \\ && \hspace{1em}=R_{-1}^{[2]}x_{\alpha_1}(-2)x_{\alpha_1}(-2) -
           x_{\alpha_2}(-2)x_{\alpha_2}(1)R_{-4}^{[1]} \nonumber \\
           &&\hspace{2em}+2x_{\alpha_2}(-2)x_{\alpha_2}(1)x_{\alpha_1}(-3)
           x_{\alpha_1}(-1) \nonumber 
\end{eqnarray} where both formulas are modulo 
         $\widetilde{U(\bar{\goth{n}}) \bar{\goth{n}}^{+}}$.  Therefore, we have
         obtained the decomposition $M=M_1+M_2+M_3$, where $M_1$ is a
         linear combination of
         $$x_{\alpha_2}(-1)x_{\alpha_2}(1)x_{\alpha_1}(-4)x_{\alpha_1}(-1)
\; \; \mbox{and} \; \;
x_{\alpha_2}(-2)x_{\alpha_2}(1)x_{\alpha_1}(-3)x_{\alpha_1}(-1),
         $$
         $M_2 \in {\mathcal J}$ and $M_3 \in
         \widetilde{U(\bar{\goth{n}})\bar{\goth{n}}^{+}}$.
 \end{example} 

 The proof of the next result is completely analogous to the proof of
 Lemma \ref{p1}:

  \begin{lemma} \label{p1-12} Let $r_1, r_2 \geq 1$ and $m_{1,2},
    \dots , m_{r_2,2}, m_{1,1}, \dots , m_{r_1,1}$ be integers such
    that $m_{1,2} \leq \cdots \leq m_{r_2,1} $ and $m_{1,1}
    \leq \cdots \leq m_{r_1,1} $, and consider
$$
    M=x_{\alpha_2}(m_{1,2})\cdots
    x_{\alpha_2}(m_{r_2,2})x_{\alpha_1}(m_{1,1}) \cdots
    x_{\alpha_1}(m_{r_1,1}).
    $$
  \begin{enumerate}
  \item If $m_{r_1,1} \leq -2$ and $m_{r_2,1} \leq r_1-1$ we have
    $$
    M=M_1+M_2+M_3+M_4,
    $$
    where $M_1$ is a linear combination of monomials which satisfy the
    difference two condition and which do not end in
    $x_{\alpha_1}(-1)$, $M_2 \in {\mathcal J}$, $M_3 \in
    \widetilde{U(\bar{\goth{n}})\bar{\goth{n}}^{+}}$ and $M_4 \in
    U(\bar{\goth{n}}) x_{\alpha_1}(-1)$. We also have that $M_1, M_2,
    M_3, M_4$ and $M$ have equal weights.

    \item
    If $m_{r_1,1} \leq -1$ and $m_{r_2, 2} \leq r_2-2$ then
     $$
    M=M_1+M_2+M_3+M_4,
    $$
    where $M_1$ is a linear combination of monomials which satisfy the
    difference two condition and which do not end in
    $x_{\alpha_2}(-1)$, $M_2 \in {\mathcal J}$, $M_3 \in
    \widetilde{U(\bar{\goth{n}})\bar{\goth{n}}^{+}}$ and $M_4 \in
    U(\bar{\goth{n}}) x_{\alpha_2}(-1)$. We also have that $M_1, M_2,
    M_3, M_4$ and $M$ have equal weights.
\end{enumerate}
In particular, (1) and (2) hold for homogeneous elements. $\; \; \;
\Box$
\end{lemma}

 The next result shows, in particular, that if the element $M_1$
 {}from Lemma \ref{p1} belongs to the ideal $I_{\Lambda_0}$ then it
 must be zero (recall (\ref{ker})).  To prove this result, we follow
 an idea developed in \cite{G1}, which involves the linear maps ${\cal
 Y}_c(e^{\lambda_i}, x)$ for $i=1,2$. See the end of Section 2 above
 for the relevant properties of these maps.

   \begin{lemma} \label{p2} Let $r_1, r_2\geq 1$ and ${\textbf
  m_1}=(m_{1,1}, \dots , m_{r_1,1})$ and ${\textbf m_2}=(m_{1,2},
  \dots m_{r_2,2})$ be sequences such that they satisfy the difference
  two condition, and $m_{r_1,1} \leq -1$ and $m_{r_2,2} \leq r_1 -1$.
  If the following homogeneous element of positive weight:
  \begin{equation} \label{indep} \sum_{\textbf m_2, \; \textbf m_1}
  A_{\textbf m_2, \; \textbf m_1} x_{\alpha_2}(m_{1,2}) \cdots
  x_{\alpha_2}(m_{r_2,2}) x_{\alpha_1}(m_{1,1}) \cdots
  x_{\alpha_1}(m_{r_1,1}) \end{equation} belongs to the ideal
  $I_{\Lambda_0}$ then all constants $A_{\textbf m_2, \; \textbf m_1}$
  equal zero.  \end{lemma} {\em Proof:} Since (\ref{indep}) belongs to
  $I_{\Lambda_0}$,
\begin{equation}  \label{12}
 \sum_{{\textbf m_2}, \; {\textbf m_1}} A_{{\textbf m_2}, \; {\textbf
  m_1}} x_{\alpha_2}(m_{1,2}) \cdots x_{\alpha_2}(m_{r_2,2})
  x_{\alpha_1}(m_{1,1}) \cdots x_{\alpha_1}(m_{r_1,1}) \cdot
  v_{\Lambda_0}=0.  \end{equation} We will prove by induction on
  weight $k= -(m_{1,2}+ \cdots m_{r_2, 2}+m_{1,1} + \cdots +
  m_{r_1,1})$ that all coefficients $A_{{\textbf m_2}, \; {\textbf
  m_1}}$ are zero. Note that $k \geq r_1^2+r_2^2-r_1r_2 >0$ for all
  $r_1, r_2 \geq 1$.  Equation (\ref{12}) is equivalent with
  \begin{eqnarray} \label{12-1} && \sum_{{\textbf m_2}, \; {\textbf
  m_1'}} A_{{\textbf m_2}, \; {\textbf m_1'}} x_{\alpha_2}(m_{1,2})
  \cdots x_{\alpha_2}(m_{r_2,2}) x_{\alpha_1}(m_{1,1}) \cdots
  x_{\alpha_1}(-1) \cdot v_{\Lambda_0} \\ && + \sum_{{\textbf m_2}, \;
  {\textbf m_1''}} A_{{\textbf m_2}, \; {\textbf m_1''}}
  x_{\alpha_2}(m_{1,2}) \cdots x_{\alpha_2}(m_{r_2,2})
  x_{\alpha_1}(m_{1,1}) \cdots x_{\alpha_1}(m_{r_1,1}) \cdot
  v_{\Lambda_0}=0, \nonumber \end{eqnarray} where ${\textbf m_1'}=
  (m_{1,1}, \dots, -1)$ and ${\textbf m_1''}= (m_{1,1}, \dots, m_{r_1,
  1})$ with $m_{r_1,1} \leq -2$. The sum of the indices of ${\textbf
  m_1'}$ equals the sum of the the indices of ${\textbf m_1''}$, and
  both equal the sum of the indices of the sequence ${\textbf m_1}$.

  We apply the linear map ${\cal Y}_c (e^{\lambda_1}, x)$ to both
  sides of (\ref{12-1}), use the properties (\ref{int1}) and
  (\ref{int2}) of this map, and then we use (\ref{prod-oper}) to move
  $e^{\lambda_1}$ to the left of the expression. Finally multiply by
  the inverse of $e^{\lambda_1}$ and obtain:
\[
    \sum_{{\textbf m_2} \; {\textbf m_1''}} A_{{\textbf m_2}, \;
    {\textbf m_1''}} x_{\alpha_2}(m_{1,2}) \cdots
    x_{\alpha_2}(m_{r_2,2}) x_{\alpha_1}(m_{1,1}+1) \cdots
    x_{\alpha_1}(m_{r_1,1}+1) \cdot v_{\Lambda_0}=0, \] where
    $m_{r_1,1} \leq -2$. Since this is a sum of homogeneous elements
    of weight less than $k$ applied to $v_{\Lambda_0}$ by induction
    and by (\ref{direct sum j}) we obtain $A_{{\textbf m_2}, \;
    {\textbf m_1''}}=0$.

As for the remaining part of (\ref{12-1}),
\begin{equation} \label{12-2}
\sum_{{\textbf m_2}, \; {\textbf m_1'}} A_{{\textbf m_2}, \; {\textbf
m_1'}} x_{\alpha_2}(m_{1,2}) \cdots x_{\alpha_2}(m_{r_2,2})
x_{\alpha_1}(m_{1,1}) \cdots x_{\alpha_1}(-1) \cdot v_{\Lambda_0}=0
\end{equation} 
we use the fact that $x_{\alpha_1}(-1) \cdot
v_{\Lambda_0}=e^{\alpha_1}$ and move $e^{\alpha_1}$ to the left of the
expression (\ref{12-2}) so that we obtain
\begin{equation} \label{12-3}
\sum_{{\textbf m_2}, \; {\textbf m_1'}} A_{{\textbf m_2}, \; {\textbf
 m_1'}} x_{\alpha_2}(m_{1,2}-1) \cdots x_{\alpha_2}(m_{r_2,2}-1)
 x_{\alpha_1}(m_{1,1}+2) \cdots x_{\alpha_1}(m_{r_1-1, 1}+2) \cdot
 v_{\Lambda_0}=0.
\end{equation}
Now we apply the linear map ${\cal Y}_c(e^{\lambda_2}, x)$ to the both
sides of (\ref{12-3}), use its properties and move $e^{\lambda_2}$ to
the left, so that we get
$$
\sum_{{\textbf m_2}, \; {\textbf m_1'}} A_{{\textbf m_2}, \; {\textbf
m_1'}} x_{\alpha_2}(m_{1,2}) \cdots x_{\alpha_2}(m_{r_2,2})
x_{\alpha_1}(m_{1,1}+2) \cdots x_{\alpha_1}(m_{r_1-1, 1}+2) \cdot
v_{\Lambda_0}=0.
$$
 The induction hypothesis combined with (\ref{direct sum j}) implies
 $A_{{\textbf m_2}, \; {\textbf m_1'}}=0$. This concludes the proof of
 the lemma. $\; \; \; \Box$ \vspace{1em}

 We have an analogous result and its proof reduces to the proof of
 that of Lemma \ref{p2}:

 \begin{lemma}\label{p2-12} Let $r_1, r_2\geq 1$ and ${\textbf
 m_1}=(m_{1,1}, \dots , m_{r_1,1})$ and ${\textbf m_2}=(m_{1,2}, \dots
 m_{r_2,2})$ be sequences such that they satisfy the difference two
 condition \begin{enumerate} \item Assume that $m_{r_1,1} \leq -2$ and
 $m_{r_2,2} \leq r_1 -1$.  If the following homogeneous element of
 positive weight: \begin{equation} \label{ind1} \sum_{\textbf m_2, \;
 \textbf m_1} A_{\textbf m_2, \; \textbf m_1} x_{\alpha_2}(m_{1,2})
 \cdots x_{\alpha_2}(m_{r_2,2}) x_{\alpha_1}(m_{1,1}) \cdots
 x_{\alpha_1}(m_{r_1,1}) \end{equation} belongs to the ideal
 $I_{\Lambda_1}$ then all constants $A_{\textbf m_2, \; \textbf m_1}$
 equal zero.

  \item Assume that $m_{r_1,1} \leq -1$ and $m_{r_2, 2} \leq r_1-2$.
  If the following homogeneous element of positive weight:
  \begin{equation} \label{ind2} \sum_{\textbf m_2, \; \textbf m_1}
  A_{\textbf m_2, \; \textbf m_1} x_{\alpha_2}(m_{1,2}) \cdots
  x_{\alpha_2}(m_{r_2,2}) x_{\alpha_1}(m_{1,1}) \cdots
  x_{\alpha_1}(m_{r_1,1}) \end{equation} belongs to the ideal
  $I_{\Lambda_2}$ then all constants $A_{\textbf m_2, \; \textbf m_1}$
  equal zero.
\end{enumerate}
  \end{lemma} {\em Proof:} Since (\ref{ind1}) lies in $I_{\Lambda_1}$
 and (\ref{ind2}) lies in $I_{\Lambda_2}$ then
\begin{equation}  
   \sum_{\textbf m_2, \; \textbf m_1} A_{\textbf m_2, \; \textbf m_1}
  x_{\alpha_2}(m_{1,2}) \cdots x_{\alpha_2}(m_{r_2,2})
  x_{\alpha_1}(m_{1,1}) \cdots x_{\alpha_1}(m_{r_1,1}) \cdot
  v_{\Lambda_1}=0 \end{equation} and \begin{equation} \sum_{\textbf
  m_2, \; \textbf m_1} A_{\textbf m_2, \; \textbf m_1}
  x_{\alpha_2}(m_{1,2}) \cdots x_{\alpha_2}(m_{r_2,2})
  x_{\alpha_1}(m_{1,1}) \cdots x_{\alpha_1}(m_{r_1,1}) \cdot
  v_{\Lambda_2}=0.  \end{equation} By using (\ref{vectors}) and
  (\ref{prod-oper}) we get: \begin{equation} \label{ind11}
  \sum_{\textbf m_2, \; \textbf m_1} A_{\textbf m_2, \; \textbf m_1}
  x_{\alpha_2}(m_{1,2}) \cdots x_{\alpha_2}(m_{r_2,2})
  x_{\alpha_1}(m_{1,1}+1) \cdots x_{\alpha_1}(m_{r_1,1}+1) \cdot
  v_{\Lambda_0}=0 \end{equation} and \begin{equation} \label{ind22}
  \sum_{\textbf m_2, \; \textbf m_1} A_{\textbf m_2, \; \textbf m_1}
  x_{\alpha_2}(m_{1,2}+1) \cdots x_{\alpha_2}(m_{r_2,2}+1)
  x_{\alpha_1}(m_{1,1}) \cdots x_{\alpha_1}(m_{r_1,1}) \cdot
  v_{\Lambda_0}=0 \end{equation} which imply that all coefficients are
  zero by using Lemma \ref{p2}.  $\Box$ \vspace{1em}

Finally, we have all the necessary prerequisites to prove the theorem
which describes the left ideals $I_{\Lambda_i}$ of
$U(\bar{\goth{n}})$ for $i=0,1,2$.  

{\em Proof of Theorem \ref{presentation sl3}}: The square of the
vertex operator $Y(e^{\alpha_j}, x)$ with $j=1, 2$ is well defined
(the components $x_{\alpha_j}(m)$, $m \in \mathbb{Z}$ of
$Y(e^{\alpha_j}, x)$ commute) and equals zero on $L(\Lambda_i)$ for
$i=0,1,2$. In particular, the square operator $Y(e^{\alpha_j}, x)^2$
is zero on $W(\Lambda_i)$. The expansion coefficients of
$Y(e^{\alpha_1}, x)^2$ and $Y(e^{\alpha_2}, x)^2$ are $R_t^{[1]}$ and
$R_t^{[2]}$ for $t \in \mathbb{Z}$:
\begin{equation}
   Y(e^{\alpha_1}, x)^2 = \sum_{t \in \mathbb{Z}} \left (
  \sum_{m_1+m_2=t} x_{\alpha_1}(m_1)x_{\alpha_1}(m_2) \right )
  x^{-t-2} \end{equation} and \begin{equation} Y(e^{\alpha_2}, x)^2=
  \sum_{t \in \mathbb{Z}} \left ( \sum_{m_1+m_2=t}
  x_{\alpha_2}(m_1)x_{\alpha_2}(m_2) \right ) x^{-t-2}.
  \end{equation} Thus
\begin{equation} \label{1}
{\mathcal J} \subset I_{\Lambda_i} \; \; \; \mbox{modulo} \; \; \; 
\widetilde{U(\bar{\goth{n}})\bar{\goth{n}}^{+}} \; \; \;  \mbox{for} \; \; \;  i=0,1,2.
\end{equation}
We also have
$$
x_{\alpha_1}(-1) \cdot v_{\Lambda_1}=0 \; \; \; \mbox{and} 
\; \;  x_{\alpha_2} (-1) \cdot v_{\Lambda_2}=0
$$
(recall (\ref{vectors}) and (\ref{prod-oper})). Hence
\begin{equation} \label{3}
U(\bar{\goth{n}})x_{\alpha_j}(-1) \subset I_{\Lambda_j} \; \;
\mbox{for} \; \; j=1,2.
\end{equation} 
Now by (\ref{1}) with $i=1,2$ and (\ref{3}) we have 
\begin{equation} {\mathcal J}+ U(\bar{\goth{n}})
   x_{\alpha_j}(-1) \subset I_{\Lambda_j} \; \; \; \mbox{modulo} \; \; \; 
\widetilde{U(\bar{\goth{n}})\bar{\goth{n}}^{+}}
\end{equation} for $j=1,2$.

 We now prove the inclusion \begin{equation} \label{inclusion2}
 I_{\Lambda_0} \subset {\mathcal J} \; \; \; 
\mbox{modulo} \; \; \;  \widetilde{U(\bar{\goth{n}}) \bar{\goth{n}}^{+}}.
\end{equation} Let us take an element $\Phi$ of $I_{\Lambda_0}$. By
Remark \ref{ideals} we may and do assume that $\Phi$ is homogeneous
with respect to the weight grading. We shall prove that $\Phi$ is 
in ${\mathcal J}$ modulo 
$\widetilde{U(\bar{\goth{n}})\bar{\goth{n}}^{+}}$. Suppose that $\mbox{wt} \;
\Phi >0$. Otherwise, by Remark \ref{negative weight} we have $\Phi
\in U(\bar{\goth{n}}) \bar{\goth{n}}^{+}$.  Lemma \ref{truncation}
applied to the homogeneous element $\Phi$ implies that
$\Phi=\Phi'+\Phi''$, with $\Phi'' \in
U(\bar{\goth{n}})\bar{\goth{n}}^{+}$ and
\begin{equation}
\Phi'= \sum_{{\textbf m_2}, \; {\textbf m_1}} A_{{\textbf m_2},
{\textbf m_1}} x_{\alpha_2}(m_{1,2}) \cdots
x_{\alpha_2}(m_{r_2,2})x_{\alpha_1}(m_{1,1}) \cdots x_{\alpha_1}
(m_{r_1,1}),
\end{equation}where the sequences 
${\textbf m_2}=(m_{1,2}, \dots , m_{r_2,2})$ and ${\textbf m_1}=
(m_{1,1}, \dots, m_{r_1,1})$ have the properties:
$$
m_{1,2} \leq \cdots \leq m_{r_2,2} \leq r_1-1\; \; \mbox{and} \; \;
m_{1,1} \leq \cdots \leq m_{r_1,1} \leq -1.
$$
Note that both $\Phi'$ and $\Phi''$ are homogeneous and have the same
weight as $\Phi$. In Lemma \ref{p1} we have shown that
$$
\Phi'= \Phi_1' + \Phi_2' + \Phi_3',
$$
where $\Phi_1'$ is a linear combination of monomials that satisfy the
difference two condition, $\Phi_2' \in {\mathcal J}$ and $\Phi_3' \in
\widetilde{U(\bar{\goth{n}})\bar{\goth{n}}^{+}} $. Moreover,
$\Phi_1'$, $\Phi_2'$ and $\Phi_3'$ have the same weight as
$\Phi'$. Since
$$
\Phi_1'=\Phi-\Phi''-\Phi_2'-\Phi_3',
$$
by using (\ref{1}) for $i=0$ we get $\Phi_1' \in I_{\Lambda_0}$ modulo
$\widetilde{U(\bar{\goth{n}})\bar{\goth{n}}^{+}}$.  Lemma \ref{p2} implies $\Phi_1'=0$. In conclusion
we obtain $\Phi \in {\mathcal J}$ modulo
$\widetilde{U(\bar{\goth{n}})\bar{\goth{n}}^{+}}$.

By using Lemmas \ref{truncation 12}, \ref{p1-12}, \ref{p2-12} and
imitating the proof of the inclusion (\ref{inclusion2}) we obtain the
inclusions
\begin{equation}
I_{\Lambda_j} \subset {\mathcal J}+
U(\bar{\goth{n}}) x_{\alpha_j}(-1) \; \; \; \mbox{modulo} \; \; \; 
\widetilde{U(\bar{\goth{n}})\bar{\goth{n}}^{+}}
\end{equation}
for $j=1,2$, and thus we prove our theorem.
$\Box$
\vspace{1em}

Assume now that $\goth{g}=\goth{sl}(l+1)$ with $l \geq 3$. As we
mentioned earlier, the proof of a presentation for the principal
subspaces $W(\Lambda_i)$ of the standard modules $L(\Lambda_i)$, where
$i=0, \dots, l$ is completely analogous to the proof of Theorem
\ref{presentation sl3}.  We observe that all the results that we have
just obtained in the case $l=2$ have obvious generalizations for all
$l \geq 2$. Indeed, the only simple root vectors that fail to commute
with one another are those corresponding to consecutive roots on the
Dynkin diagram (recall that $[x_{\alpha_i}(m), x_{\alpha_j}(n)]=0$ if
$|i-j| >1$ for $i, j =1, \dots, l$ and $m, n \in \mathbb{Z}$). Thus we
have results similar to Lemmas \ref{truncation} and \ref{truncation
  12}. Furthermore, by using the expressions $R_t^{[j]}$ and
$R_{t;m}^{[j]}$ for $1 \leq j \leq l$ and $m, t \in \mathbb{Z}$ we
obtain results analogous to Lemmas \ref{p1} and \ref{p1-12}. The
linear maps ${\cal Y}_c(e^{\lambda_j}, x)$, $j=1, \dots, l$ and their
relevant properties ((\ref{int1}) and (\ref{int2})) are the main
ingredients in order to prove results similar to Lemmas \ref{p2} and
\ref{p2-12}. Recall {}from the beginning of this section the two-sided
ideal ${\mathcal J}$ of $\widetilde{U(\bar{\goth{n}})}$ generated by $R^{[j]}_t$ for $j=1, \dots, l$ and $t \in \mathbb{Z}$.  Thus we
obtain a presentation of the principal subspaces of the level $1$
standard $\widehat{\goth{sl}(l+1)}$-modules:

\begin{theorem} \label{presentation}
The annihilator of the highest weight vector of $L(\Lambda_0)$ in
 $U(\bar{\goth{n}})$ is described by: \begin{equation}
 I_{\Lambda_0}\equiv {\mathcal J} \; \; \; \mbox{modulo} \; \; \; 
\widetilde{U(\bar{\goth{n}})\bar{\goth{n}}^{+}}.  \end{equation}
 The annihilators of the highest weight vectors of $L(\Lambda_j)$ in
 $U(\bar{\goth{n}})$ are described as follows: \begin{equation}
 I_{\Lambda_j} \equiv {\mathcal J}+
U(\bar{\goth{n}})x_{\alpha_j}(-1) 
\; \; \; \mbox{modulo} \; \; \; 
\widetilde{U(\bar{\goth{n}})\bar{\goth{n}}^{+}}
\end{equation} for any $1 \leq j
 \leq l$.  $\; \; \; \Box$ \end{theorem}

A crucial consequence of Theorem \ref{presentation} is the following
corollary, which gives the discrepancy among the ideals
$I_{\Lambda_i}$ for $i=0, \dots, l$. This result will be used in the
proof of the exactness of our sequences of maps among principal
subspaces in the next section.

 \begin{corollary} \label{discrepancy} For any integer $j$ with $1
\leq j \leq l$ we have
\begin{equation} \label{important formula} 
I_{\Lambda_j} = I_{\Lambda_0} + U(\bar{\goth{n}})x_{\alpha_j}(-1). \; \; \; 
\Box
\end{equation}
\end{corollary} 

\begin{remark}
\rm All the results proved in this section simplify substantially in
the case of $\widehat{\goth{sl}(2)}$, i.e., $l=1$. By combining these
results we obtain a proof of a presentation of the principal subspaces
of the level $1$ standard $\widehat{\goth{sl}(2)}$-modules that is
different than the one given in \cite{FS1}-\cite{FS2}.
\end{remark}

 \section{Exact sequences and recursions}
 \setcounter{equation}{0}

In this section we obtain a complete set of recursions ($q$-difference
equations) that characterize the graded dimensions of the principal
subspaces of the level one standard modules for $\widehat
{\goth{sl}(l+1)}$ by setting up $l$ exact sequences for these
principal subspaces.

The space $V_P$ has certain gradings. The natural action of the
Virasoro algebra operator $L(0)$ on $V_P$ gives a grading by {\it
weight} and this is a $ \mathbb{Q}$-grading. We have $ \mbox{wt}(e^{
\lambda})= \frac{1}{2} \langle \lambda, \lambda \rangle \in
\mathbb{Q}$ for any $ \lambda \in P$. The weight of $h(-n)$, viewed as
either an operator or as element of $U(\widehat{\goth{h}})$, equals
$n$ for any $n \in \mathbb{Z}$ and $h \in \goth{h}$. See Section 3 for
further details about the grading by weight.  The space $V_P$ has also
gradings by {\it charge}, given by the eigenvalues of the operators $
\lambda_{j}$ (thought as $\lambda_{j}(0)$) for $j=1, \dots, l$. These
gradings are $\mathbb{Q}$-gradings compatible with the grading by
weight. The charge with respect to $\lambda_j$ of $e^{\lambda_k}$,
viewed as an operator or as an element of $L(\Lambda_k)$, is $\langle
\lambda_j, \lambda_k \rangle$.  We shall restrict these gradings to
the principal subspaces $W(\Lambda_{i})$ for $i=0, \dots, l$.

Let $r_1, \dots, r_l$ be nonnegative integers. For $m_{1, l}, \dots,
m_{r_l, l}, \dots, m_{1,1}, \dots, m_{r_1, 1}$ integers the elements:
\begin{equation}
x_{\alpha_l}(m_{1,l}) \cdots x_{\alpha_l}(m_{r_l,l}) \cdots
x_{\alpha_1}(m_{1,1}) \cdots x_{\alpha_1}(m_{r_1,1}) \cdot
v_{\Lambda_0} \in W(\Lambda_0)
\end{equation}
and
\begin{equation}
x_{\alpha_l}(m_{1,l}) \cdots x_{\alpha_l}(m_{r_l,l}) \cdots
x_{\alpha_1}(m_{1,1}) \cdots x_{\alpha_1}(m_{r_1,1}) \cdot
v_{\Lambda_k} \in W(\Lambda_k)
\end{equation} 
have weights $-m_{1,l} - \cdots-m_{r_1, 1}$ and $-m_{1,l}-\cdots
-m_{r_1,1}+\frac{1}{2} \langle \lambda_k, \lambda_k \rangle$,
respectively. Their charges with respect to the operator $\lambda_j$
are $r_j$ and $r_j+ \langle \lambda_k, \lambda_j \rangle$ for all $j,k
=1, \dots, l$.

Now for any $i=0, \dots, l$ we consider the graded dimensions of the
principal subspaces $W(\Lambda_i)$ (the generating functions of the
dimensions of the homogeneous subspaces of $W(\Lambda_i)$):
\begin{equation}
\chi _i(x_1, \dots, x_l; q)= \mbox{dim}_{*} (W(\Lambda_i), x_1, \dots
, x_l;q)= tr|_{W(\Lambda_i)}x_{1}^{ \lambda_ 1} \dots x_{l}^{
\lambda_l}q^{L(0)},
\end{equation}
where $x_1, \dots , x_l$ and $q$ are formal variables. 
These generating functions are also called the characters
 of the principal subspaces $W(\Lambda_i)$. 

 Notice that 
 $$
 \chi_{0}(x_1, \dots, x_l;q) \in \mathbb{C}[[x_1, \dots, x_l;q]].
 $$
To avoid the multiplicative factors $x_1^{\langle \lambda_1, \lambda_j
 \rangle} \cdots x_l^{\langle \lambda_l, \lambda_j \rangle } q^{
 \frac{1}{2} \langle \lambda_j, \lambda_j \rangle}$, where $j=1,
 \dots, l$ we use slightly modified graded dimensions as follows:
 \begin{equation} \label{chi'} \chi_{j}'(x_1, \dots, x_l; q)=
 x_1^{-\langle \lambda_1, \lambda_j \rangle} \cdots x_l^{-\langle
 \lambda_l, \lambda_j \rangle } q^{- \frac{1}{2} \langle \lambda_j,
 \lambda_j \rangle} \chi_j (x_1, \dots, x_l; q).  \end{equation} Thus
 we have
$$
\chi_j'(x_1, \dots, x_l;q) \in \mathbb{C}[[x_1, \dots, x_l,q]].
$$ 
We shall also use the following notation for the homogeneous subspaces
of $W(\Lambda_j)$ with respect to the above gradings:
\begin{equation} \label{W'}
W(\Lambda_j)'_{r_1, \dots, r_l; k}=W(\Lambda_j)_{r_1+\langle
\lambda_1, \lambda_j \rangle, \dots, r_l+ \langle \lambda_l, \lambda_j
\rangle ; k+ 1/2 \langle \lambda_j, \lambda_j \rangle},
\end{equation}
where $j=1, \dots , l$ and $r_1, \dots, r_l, k \geq 0$.

We now use certain maps corresponding to the intertwining operators
${\cal Y}(e^{\lambda_j}, x)$, where $j=1, \dots, l$ in order to get
relations between the graded dimensions of $W(\Lambda_j)$ and
$W(\Lambda_0)$.

Consider the linear maps
 $$
 e^{\lambda_j}: V_P \longrightarrow V_P
 $$
 for any $j=1, \dots, l$.
 Clearly, they are isomorphisms with $e^{-\lambda_j}$ as their inverses.

\begin{proposition} \label{characters}
The following relations among the graded dimensions of the principal
subspaces $W( \Lambda_0)$, $W(\Lambda _1), \dots , W(\Lambda_l)$ hold:
\begin{equation}  \label{relation}
 \chi_1 '(x_1, x_2, \dots, x_l ; q)=  \chi _0 (x_1q, x_2, \dots, x_l; q), 
 \end{equation}
 \begin{equation} \label{second relation}
  \chi _2'( x_1, x_2, \dots, x_l ; q) = \chi _0 (x_1, x_2q, \dots, x_l; q),
  \end{equation}
  $$
  \vdots
  $$
  \begin{equation} \label{third relation} 
\chi _l'( x_1, x_2, \dots,
x_l; q) = \chi _0(x_1, x_2, \dots, x_lq;q).
\end{equation}
\end{proposition}
{ \em Proof:} We restrict the isomorphism $e^{\lambda_1}$ to
the principal subspace $W(\Lambda_0)$ of $L(\Lambda_0)$.
We have
$$
e^{\lambda_1} x_{\alpha_1}(m)= x_{\alpha_1}(m-1) e^{\lambda_1},
$$
$$
e^{\lambda_1} x_{\alpha_j}(m)=x_{\alpha_j}(m) e^{\lambda_1}, \; \;
\mbox{for} \; \; j=1, \dots, l
$$
on $V_P$ for $m \in \mathbb{Z}$. We also have
$$
e^{\lambda_1} \cdot v_{\Lambda_0}= v_{\Lambda_1}.
$$
 Thus \begin{eqnarray} &&e^{\lambda_1} (x_{\alpha_l}(m_{1,l}) \cdots
 x_{\alpha_l}(m_{r_l,l}) \cdots x_{\alpha_1}(m_{1,1}) \cdots
 x_{\alpha_1}(m_{r_1,1}) \cdot v_{\Lambda_0}) \\ &&=
 x_{\alpha_l}(m_{1,l}) \cdots x_{\alpha_l}(m_{r_l,l}) \cdots
 x_{\alpha_1}(m_{1,1}-1) \cdots x_{\alpha_1}(m_{r_1,1}-1) \cdot
 v_{\Lambda_1} \nonumber \end{eqnarray} for $m_{1,l}, \dots ,
 m_{r_l,l}, \dots, m_{1,1}, \dots, m_{r_1,1} \in \mathbb{Z}$. Then we
 obtain
\begin{equation} \label{e-1}
e^{ \lambda_1} : W(\Lambda_0)  \longrightarrow  W(\Lambda_1), 
\end{equation}
a linear isomorphism.

The map $e^{\lambda_1}$ does not preserve weight and charges.  Let
$W(\Lambda_0)_{r_1, \dots, r_l; k}$ be an homogeneous subspace of
$W(\Lambda_0)$ with $r_1, \dots, r_l, k \geq 0$. We observe that the
map $e^{\lambda_1}$ increases the charge corresponding to $\lambda_j$
by $\langle \lambda_1, \lambda_j \rangle $ for all $j=1, \dots, l$.
Now for any $w \in W(\Lambda_0)_{r_1, \dots, r_l; k}$ the element
$e^{\lambda_1}(w) $ has weight $k+r_1+ 1/2 \langle \lambda_1,
\lambda_1 \rangle$. Thus we obtain the following isomorphism of
homogeneous spaces: \begin{equation} e^{ \lambda_1} :
W(\Lambda_0)_{r_1, \dots , r_l;k } \longrightarrow W(\Lambda_1)_{ r_1,
\dots , r_l; k+r_1}' \end{equation} (recall (\ref{W'}) for notation),
which proves the relation between the graded dimensions of
$W(\Lambda_1)$ and $W(\Lambda_0)$
$$ \chi_1' (x_1, x_2, \dots x_l; q)= \chi _0 (x_1q, x_2, \dots, x_l;
q).$$

 The remaining relations of this proposition are proved similarly by
 using the restrictions of the linear isomorphisms $e^{\lambda_j}$ to
 $W(\Lambda_0)$ for all $j=2, \dots, l$.  $\; \; \; \Box$ \vspace{1em}

 \begin{remark} \rm One can view the map $e^{\lambda_j}:W(\Lambda_0)
 \longrightarrow W(\Lambda_j)$ {}from the previous proposition as
 essentially the ``constant factor'' of the intertwining operator
 ${\cal Y}(e^{\lambda_j}, x)$ for any $j=1, \dots, l$.  \end{remark}

Now we consider the weights $ \lambda^j = \alpha_j-\lambda_j$ 
 and the linear isomorphisms 
 \begin{equation} \label{j}
 e^{\lambda^j}: V_P \longrightarrow V_P
 \end{equation}
 for all $j=1, \dots, l$.
 The restriction of (\ref{j}) to $W(\Lambda_j)$ is the linear map 
 \begin{equation} \label{first}
 e^{\lambda^j}: W(\Lambda_j) \longrightarrow W(\Lambda_0).
 \end{equation}
By using (\ref{prod-oper}) and   
\begin{equation}
e^{\alpha_i} =x_{\alpha_i}(-1) \cdot
v_{\Lambda_0},
\end{equation}
which follows {}from the creation property of the vertex operator
$Y(e^{\alpha_i}, x)$, we have
\begin{eqnarray} \label{ex}
&& e^{\lambda^j} (x_{\alpha_i} (m_{1,i}) \cdots x_{\alpha_i} (m_{r_i,
i}) \cdot v_{\Lambda_j}) \\ &&=x_{\alpha_i} (m_{1,i}- \langle
\lambda^j, \alpha_i \rangle ) \cdots x_{\alpha_i} (m_{r_i, i}- \langle
\lambda^j, \alpha_i \rangle ) x_{\alpha_j}(-1)\cdot v_{\Lambda_0}
\nonumber
\end{eqnarray}
for any $i, j=1, \dots, l$ and $m_{1,i}, \dots, m_{r_i,i} \in \mathbb{Z}$.

Recall {}from the end of Section 2 the constant terms of the
 intertwining operators ${\cal Y} (e^{\lambda_j}, x) : L(\Lambda_0)
 \longrightarrow L(\Lambda_j)\{x \}$, denoted by ${\cal
 Y}_c(e^{\lambda_j}, x)$. The restriction of ${\cal
 Y}_c(e^{\lambda_j}, x)$ to $W(\Lambda_0)$ is a well-defined linear
 map between principal subspaces: \begin{equation} \label{second}
 {\cal Y}_c( e^ { \lambda_j}, x): W(\Lambda_0) \longrightarrow W(
 \Lambda_j).  \end{equation}

The main result of this paper is the following theorem, giving $l$
exact sequences of maps among principal subspaces. The relevant
properties of intertwining operators and the result about the defining
ideals for the principal subspaces {}from Section 3 (Theorem
\ref{presentation}) will be the main tools for proving this theorem.

\begin{theorem} \label{main}
Recall the linear maps $e^{\lambda^j}$ and ${\cal Y}_c(e^{\lambda_j},
x)$ introduced above (see (\ref{first}) and (\ref{second})).  The
sequences of maps between principal subspaces \begin{equation}
\label{seq1} 0 \longrightarrow W( \Lambda_1) \stackrel{e^{\lambda^1}}
\longrightarrow W( \Lambda_0) \stackrel{{\cal Y}_c (e^{\lambda_1}, x)}
\longrightarrow W( \Lambda_1) \longrightarrow 0, \end{equation}
\begin{equation} 0 \longrightarrow W(\Lambda_2)
\stackrel{e^{\lambda^2}} \longrightarrow W( \Lambda_0) \stackrel{{\cal
Y}_c (e^{\lambda_2}, x)} \longrightarrow W( \Lambda_2) \longrightarrow
0, \end{equation}
 $$
 \vdots
 $$
  \begin{equation} 0 \longrightarrow W( \Lambda_l)
 \stackrel{e^{\lambda^l}} \longrightarrow W( \Lambda_0)
 \stackrel{{\cal Y}_c (e^{\lambda_l}, x)} \longrightarrow W(
 \Lambda_l) \longrightarrow 0 \end{equation} are exact.  \end{theorem}
 {\em Proof:} We will prove that the sequence (\ref{seq1}) is
 exact. One can show similarly that the other sequences are exact.

 The map $e^{\lambda^1}$ is clearly injective. Indeed, for any $w \in
 W(\Lambda_1)$ such that $e^{\lambda^1} (w)=0$ we get $w=0$ by the
 bijectivity of $e^{\lambda^1}: V_P \longrightarrow V_P$.  Using the
 properties (\ref{int1}) and (\ref{int2}) of the linear map ${\cal
 Y}_c(e^{\lambda_1}, x)$ we obtain
\begin{eqnarray} \nonumber
&&W(\Lambda_1)=U(\bar{\goth{n}}) \cdot v_{\Lambda_1}=
U(\bar{\goth{n}}) \; {\cal Y}_c(e^{\lambda_1}, x) \cdot v_{\Lambda_0}
\nonumber \\ &&={\cal Y}_c(e^{\lambda_1},x) \; U(\bar{\goth{n}}) \cdot
v_{\Lambda_0}= {\cal Y}_c(e^{\lambda_1}, x) \; W(\Lambda_0), \nonumber
\end{eqnarray}
and this proves that ${\cal Y}_c (e^{\lambda_1}, x)$ is a surjection.

Now we show that 
\begin{equation} \label{chain}
\mbox{Im} \; e^{\lambda^1} \subset \mbox{Ker} \; {\cal Y}_c
(e^{\lambda_1},x).
\end{equation}
Let $w \in \mbox{Im} \; e^{\lambda^1}$. By (\ref{ex}) we observe that
$w=v \cdot v_{\Lambda_0}$ with $v \in
U(\bar{\goth{n}})x_{\alpha_1}(-1)$.  We apply ${\cal
Y}_c(e^{\lambda_1}, x)$ to $w$, use its properties together with
$$
x_{\alpha_1}(-1)\cdot v_{\Lambda_1}=x_{\alpha}(-1) \cdot
e^{\lambda_1}= e^{\lambda_1}x_{\alpha_1}(0)\cdot v_{\Lambda_0}=0,
$$
and thus we obtain
$$
{\cal Y}_c(e^{\lambda_1}, x) (w)=0.
$$ 
This proves the inclusion (\ref{chain}).

 It remains to prove 
\begin{equation} \label{the second inclusion}
\mbox{Ker} \; {\cal Y}_c(e^{\lambda_1}, x) \subset \mbox{Im} \;
e^{\lambda^1}.
\end{equation} 
 In order to prove this inclusion we will first characterize the
 spaces $\mbox{Ker} \; {\cal Y}_c (e^{\lambda_1}, x)$ and $\mbox{Im}
 \; e^{\lambda^1}$. Let $w$ be an element of $\mbox{Ker} \; {\cal Y}_c
 (e^{\Lambda_1}, x)$. Since $w \in W(\Lambda_0)$ we have
 $w=f_{\Lambda_0}(u)$ with $u \in U(\bar{\goth{n}})$ (recall (\ref{f
 maps})).  By using (\ref{int1}), (\ref{int2}) and (\ref{ker}) we
 obtain
$$
{\cal Y}_c (e^{\lambda_1}, x) \; f_{\Lambda_0}(u)=0
\Longleftrightarrow f_{\Lambda_1}(u)=0 \Longleftrightarrow u \in
I_{\Lambda_1}.
$$
We have just shown that the space $\mbox{Ker} \; {\cal
Y}_c(e^{\lambda_1}, x)$ is characterized as follows:
\begin{equation} \label {incl}
w=f_{\Lambda_0}(u) \in \mbox{Ker} \; {\cal Y}_c (e^{\lambda_1}, x)
\Longleftrightarrow u \in I_{\Lambda_1}.
\end{equation}  

Let $w \in \mbox{Im} \; e^{\lambda^1}$. Then $w=vx_{\alpha_1}(-1)
\cdot v_{\Lambda_0}$ with $v \in U(\bar{\goth{n}})$ (see
(\ref{ex})). On the other hand, since $w \in W(\Lambda_0)$, we have
$w=f_{\Lambda_0}(u)$ with $u \in U(\bar{\goth{n}})$ (recall (\ref{f
maps})). Thus 
$$
f_{\Lambda_0} (u)= v x_{\alpha_1}(-1) \cdot
v_{\Lambda_0}=f_{\Lambda_0}(vx_{\alpha_1}(-1)),
$$
which is equivalent with 
$$
u-vx_{\alpha_1}(-1)
\in I_{\Lambda_0}.
$$
Therefore we have obtained 
\begin{equation} \label{inclu}
w=f_{\Lambda_0}(u) \in \mbox{Im} \; e^{\lambda^1} \Longleftrightarrow
u \in U(\bar{\goth{n}})x_{\alpha_1}(-1) + I_{\Lambda_0}.
\end{equation}

The descriptions of the vector spaces $\mbox{Ker} \; {\cal
Y}_c(e^{\lambda_1}, x)$ and $\mbox{Im} \; e^{\lambda^1}$ ((\ref{incl})
and (\ref{inclu})) combined with formula (\ref{important formula})
from Corollary \ref{discrepancy} prove the inclusion (\ref{the second
inclusion}), and therefore the exactness of sequence (\ref{seq1}).
$\; \; \; \Box$
\vspace{1em}

\begin{remark} 
\rm Notice that (\ref{incl}), (\ref{inclu}) and Corollary
\ref{discrepancy} prove also the inclusion $\mbox{Im} \; e^{\lambda^1}
\subset \mbox{Ker} \; {\cal Y}_c(e^{\lambda_1}, x)$. This inclusion
follows easily, for a different reason, as we have seen in the proof
of Theorem \ref{main}.  We do not have an elementary proof of the
inclusion $\mbox{Ker} \; {\cal Y}_c(e^{\lambda_1},x) \subset \mbox{Im}
\; e^{\lambda^1}$ that does not use the description of the ideals
$I_{\Lambda_i}$ for $i=0, \dots , l$.
\end{remark}

 As the main consequence of the previous theorem we obtain a complete
 set of recursions for the graded dimension $\chi_0(x_1, \dots, x_l
 ;q)$ of $W(\Lambda_0)$:

  \begin{theorem} \label{recursions} The graded dimension of the
   principal subspace $W(\Lambda_0)$ satisfies the following
   recursions:
   $$
  \chi_{0}(x_1, \dots, x_l;q) =\chi_{0}(x_1q, x_2, x_3, \dots , x_l;
  q)+x_1q\chi_{0}(x_1q^2, x_2q^{-1}, x_3, \dots, x_{l-1}, x_l;q),
  $$
  $$
  \chi_{0}(x_1, \dots, x_l; q) = \chi_{0}(x_1, x_2q, x_3, \dots , x_l;
  q)+x_2q \chi_{0}(x_1q^{-1}, x_2q^2, x_3q^{-1}, \dots, x_l;q),
  $$
  $$
  \vdots
 $$
 $$
  \chi_0(x_1, \dots, x_l;q) = \chi_0(x_1, x_2, x_3, \dots, x_lq; q) +
  x_lq \chi_0 (x_1, x_2, x_3, \dots, x_{l-1}q^{-1}, x_lq^2;q).
  $$
   \end{theorem} 
{\em Proof:} We will prove the first recursion of
this theorem. The other recursions follow similarly.

Consider the linear map (\ref{first}) when $j=1$. Since
\begin{eqnarray} 
&& e^{\lambda^1} (x_{\alpha_1} (m_{1,1}) \cdots x_{\alpha_1} (m_{r_1,
1}) \cdot v_{\Lambda_1}) \\ &&=x_{\alpha_1} (m_{1,1}- 1) \cdots
x_{\alpha_1} (m_{r_1, 1}- 1 ) x_{\alpha_1}(-1)\cdot v_{\Lambda_0},
\nonumber
\end{eqnarray}
\begin{eqnarray} 
&& e^{\lambda^1} (x_{\alpha_1} (m_{1,2}) \cdots x_{\alpha_2} (m_{r_2,
2}) \cdot v_{\Lambda_1}) \\ &&=x_{\alpha_1} (m_{1,2}+1) \cdots
x_{\alpha_1} (m_{r_1, 2}+1 ) x_{\alpha_1}(-1)\cdot v_{\Lambda_0}
\nonumber
\end{eqnarray}
and for any $j=3, \dots, l$,
\begin{eqnarray} 
&& e^{\lambda^1} (x_{\alpha_j} (m_{1,j}) \cdots x_{\alpha_j} (m_{r_j,
j}) \cdot v_{\Lambda_1}) \\ &&=x_{\alpha_j} (m_{1,j}) \cdots
x_{\alpha_1} (m_{r_j, j}) x_{\alpha_1}(-1)\cdot v_{\Lambda_0},
\nonumber
\end{eqnarray}
where $m_{1, i}, \dots, m_{r_i,i} \in \mathbb{Z}$ and $i=1, \dots, l$,
  we have the following map between homogeneous spaces with respect to
  the weight and charge gradings: \begin{equation} \label{homog1}
  e^{\lambda^1} : W(\Lambda_1)'_{r_1, r_2, \dots, r_l;k}
  \longrightarrow W(\Lambda_0)_{r_1+1,r_2; \dots, r_l;k+r_1-r_2+1}
  \end{equation} (recall our notation (\ref{W'})).

By using (\ref{int1}) and (\ref{int2}) we observe that ${\cal
 Y}_c(e^{\lambda_1}, x)$ is of weight $1/2\langle \lambda_1, \lambda_1
 \rangle$, it increases the charge corresponding to the operator
 $\lambda_j$ by $\langle \lambda_j, \lambda_1 \rangle $ for any $j=1,
 \dots l$.  This gives us the linear map between the homogeneous
 spaces
\begin{equation}  \label{sec}
{\cal Y}_c (e^{ \lambda_1}, x): W(\Lambda_0)_{r_1, r_2, \dots, r_l; k}
\longrightarrow W(\Lambda_1)'_{r_1, r_2, \dots ,r_l; k}.
\end{equation} 

Now the exactness of the sequence (\ref{seq1}) implies 
\begin{eqnarray} \nonumber
&& \chi_{0} (x_1, x_2, x_3, \dots, x_l;q) \nonumber \\ && \hspace{1em}
= \chi_{1}'(x_1, x_2, x_3, \dots, x_l;q) \nonumber \\ && \hspace{3em}
+ x_1q\chi_1'(x_1q, x_2q^{-1}, x_3, \dots, x_l;q), \nonumber
\end{eqnarray}
and this combined with formula (\ref{relation}) gives us the
recursion
\begin{eqnarray} \nonumber
  && \chi_{0}(x_1, x_2, x_3, \dots, x_l;q) \nonumber \\ &&
   \hspace{1em} =\chi_{0}(x_1q, x_2, x_3, \dots , x_l; q) \nonumber \\
   && \hspace{3em} +x_1q\chi_{0}(x_1q^2, x_2q^{-1}, x_3, \dots,
   x_l;q).  \nonumber\; \; \Box \end{eqnarray} \vspace{1em}

 The next remark is a reformulation of the recursions in terms of the
 entries of the generalized Cartan matrix of $\goth{sl}(l+1)$.

 \begin{remark} \rm The recursions {}from Theorem \ref{recursions} can
 be expressed as follows: \begin{eqnarray} \chi_{0}(x_1, \dots, x_l;q)
 & = & \chi_{0}(x_1, \dots , (x_jq)^{\frac{a_{jj}}{2}}, \dots, x_l; q)
 \\ &+ &(x_jq)^{\frac{a_{jj}}{2}}\chi_{0}(x_1q^{a_{j1}},
 x_2q^{a_{j2}}, x_3q^{a_{j3}}, \dots, x_lq^{a_{jl}};q), \nonumber
 \end{eqnarray} where $A=(a_{ij})_{1 \leq i, j \leq l}$ is the Cartan
 matrix of the Lie algebra $\goth{sl}(l+1)$.  \end{remark}

 For any nonnegative integer $m$ we use the notation
 $$(q)_m= (1-q)(1-q^2) \cdots (1-q^m);$$
 in particular,
 $$(q)_0=1.$$

 Now by solving the recursions {}from Theorem \ref{recursions} and by
 using Proposition \ref{characters} we obtain the graded dimensions of
 $W(\Lambda_0), W(\Lambda_1), \dots, W(\Lambda_l)$:

 \begin{corollary} The graded dimensions of the principal subspaces of
  the level $1$ standard $\widehat{\goth{g}}$-modules are given by the
  following formulas: \begin{equation} \chi_0 (x_1, \dots , x_l;q)=
  \sum_{r_1, \dots, r_l \geq 0} \frac{q^{r_1^2+ \cdots + r_l^2-r_2r_1-
  \cdots -r_lr_{l-1}}}{(q)_{r_1} \cdots (q)_{r_l}} x_1^{r_1} \cdots
  x_l^{r_l} \end{equation} and \begin{equation} \label{character}
  \chi_j '(x_1, \dots , x_l;q)= \sum_{r_1, \dots, r_l \geq 0}
  \frac{q^{r_1^2+ \cdots + r_l^2+ r_j-r_2r_1- \cdots
  -r_lr_{l-1}}}{(q)_{r_1} \cdots (q)_{r_l}} x_1^{r_1} \cdots x_l^{r_l}
  \end{equation} for any $j=1, \dots, l$. In particular, we have
\begin{equation} \label{character'}
 \chi_j(x_1, \dots, x_l;q)= x_1^{\langle \lambda_1, \lambda_j \rangle}
 \cdots x_l^{\langle \lambda_l, \lambda_j \rangle} q^{\frac{1}{2}
 \langle \lambda_j, \lambda_j \rangle} \chi_j' (x_1, \dots, x_l;q).
\end{equation}

\end{corollary}
 {\em Proof:} Define $A_{r_1,r_2, \dots, r_l}(q) \in \mathbb{C}[[q]]$ by
 $$
 \chi_0(x_1, x_2, \dots , x_l;q)=\sum_{r_1, r_2, \dots r_l \geq
 0}A_{r_1,r_2, \dots, r_l}(q)x_1^{r_1}x_2^{r_2} \cdots x_l^{r_l}.
 $$
 By solving the first recursion of Theorem \ref{recursions} we obtain
 \begin{equation} \label{sequence1} A_{r_1,r_2, \dots , r_l}(q)=
 \frac{q^{2r_1-1-r_2}}{1-q^{r_1}} A_{r_1-1, r_2, \dots, r_l}(q) =
 \cdots =\frac{q^{r_1^2-r_2r_1}}{(q)_{r_1}}A_{0,r_2, \dots, r_l}(q)
\end{equation}
for all $r_1 \geq 1$ and $r_2, \dots, r_l \geq 0$.

The second recursion of Theorem \ref{recursions} gives
\begin{equation} \label{sequence2}
A_{r_1,r_2, \dots, r_l}(q)=\frac{q^{2r_2-1-(r_1+r_3)}}{1-q^{r_2}}
A_{r_1, r_2-1, \dots, r_l}(q) = \cdots =
\frac{q^{r_2^2-r_2r_1-r_3r_2}}{(q)_{r_2}} A_{r_1,0, r_3, \dots,
r_l}(q)
\end{equation}
for all $r_2 \geq 1$ and $r_1, r_3, \dots , r_l \geq 0$.
In particular, if we take $r_1=0$ in (\ref{sequence2}) we get
\begin{equation} \label{sequence3}
A_{0,r_2, \dots , r_l}(q)=\frac{q^{r_2^2-r_3r_2}}{(q)_{r_2}}A_{0,0,
r_3, \dots, r_l}(q).
\end{equation}
The formula for $A_{0,0,r_3, \dots, r_l}(q)$ will be derived {}from
the third recursion as follows:
\begin{equation}
A_{0,0, r_3, \dots , r_l}(q)=\frac{q^{r_3^2-r_4r_3}}{(q)_{r_3}}A_{0,0,
0, r_4, \dots, r_l}(q).
\end{equation}

We continue this procedure for the remaining recursions. 
Finally the last recursion of the Theorem \ref{recursions} yields
\begin{equation} \label{sequence4}
A_{r_1, \dots, r_{l-1}, r_l}(q)= \frac{q^{2r_l-1-r_{l-1}}}{1-q^{r_l}}
A_{r_1, \dots, r_{l-1}, r_l} (q)= \cdots
=\frac{q^{r_l^2-r_lr_{l-1}}}{(q)_{r_l}} A_{r_1, \dots, r_{l-1},0}(q)
\end{equation}
for all $r_1, \dots, r_{l-1} \geq 0$ and $r_l \geq 1$. Thus
\begin{equation}
A_{0, \dots ,0,r_l}(q)= \frac{q^{r_l^2}}{(q)_{r_l}} A_{0, \dots, 0}(q).
\end{equation}

Note that $A_{0,0, \dots, 0}(q)$ is the dimension of the homogeneous
space consisting of the elements of charge zero with respect to
$\lambda_1, \lambda_2, \dots, \lambda_l$. This space is
$\mathbb{C}$ and so $A_{0, \dots, 0}(q)=1$.

Therefore, by combining the above results we obtain the graded
 dimension of the principal subspace $W(\Lambda_0)$, \begin{equation}
 \label{graded-0} \chi_0 (x_1, \dots , x_l;q)= \sum_{r_1, \dots, r_l
 \geq 0} \frac{q^{r_1^2+ \cdots + r_l^2-r_2r_1- \cdots
 r_lr_{l-1}}}{(q)_{r_1} \cdots (q)_{r_l}} x_1^{r_1} \cdots x_l^{r_l}.
 \end{equation} Now the formulas (\ref{character}) and
 (\ref{character'}) follow {}from Proposition \ref{characters} and
 {}from the equation (\ref{chi'}), respectively. $\; \; \; \Box$

\vspace{.4in}

\noindent {\small \sc Department of Mathematics, Rutgers University,
Piscataway, NJ 08854}

\vspace{.2in}

\noindent Current address:\\
\noindent{\small \sc Department of Mathematics, Ohio State University,
Columbus, OH 43210} \\
{\em E--mail address}: calinescu@math.ohio-state.edu

\end{document}